\documentclass[11pt]{amsart}

\usepackage{mathrsfs}
\usepackage{amsmath, amscd, amsthm, amssymb, amsfonts, verbatim,subfigure}

\usepackage[all]{xy}
\usepackage{graphicx}

\title[TCFTs and gauge theories]{Topological conformal field theories and gauge theories}
\author{Kevin Costello}
\address{Department of Mathematics \\ University of Chicago}
\email{costello@math.uchicago.edu}

\date{}

\newcommand{\dens}[1]{\op{Densities} (#1)}
\newcommand{\A}{\mc A}

\newcommand{\tr}{\triangle}

\newcommand{\til}{\widetilde}
\newcommand{\mscr}{\mathscr}
\renewcommand{\det}{\operatorname{det}}

\newcommand{\Comp}{\op{Comp}_{\C}}

\newcommand{\br}{\overline}

\newcommand{\iso}{\cong}
\newcommand{\C}{\mathbb C}

\newcommand{\Z}{\mathbb Z}
\newcommand{\defeq}{\overset{\text{def}}{=}}
\newcommand{\into}{\hookrightarrow}

\newcommand{\op}{\operatorname}

\newcommand{\mbb}{\mathbb}

\newcommand{\mc}{\mathcal}

\newcommand{\ip}[1]{\left\langle #1 \right\rangle}
\newcommand{\abs}[1]{\left| #1 \right|}

\newcommand{\R}{\mbb R}
\renewcommand{\d}{\mathrm{d}}
\renewcommand{\epsilon}{\varepsilon}

\renewcommand{\Im}{\op{Im}}
\newcommand{\YM}{\op{YM}}
\newcommand{\xto}{\xrightarrow}

\DeclareMathOperator{\Aut}{Aut} \DeclareMathOperator{\End}{End}
 
 \DeclareMathOperator{\Hom}{Hom}

\DeclareMathOperator{\Met}{Met} \DeclareMathOperator{\Vol}{Vol}
\DeclareMathOperator{\Tr}{Tr} 
\DeclareMathOperator{\Or}{Or}\DeclareMathOperator{\Ker}{Ker}
\DeclareMathOperator{\Mat}{Mat}

\newtheorem{theorem}{Theorem}[subsection]
\newtheorem{thm-def}{Theorem/Definition}[theorem]
\newtheorem{proposition}[theorem]{Proposition}

\newtheorem{lemma}[theorem]{Lemma}

\newtheorem{definition}[theorem]{Definition}

\numberwithin{equation}{subsection}

\newtheoremstyle{rem}
  {4pt}
  {5pt}
  {}
  {}
  {\itshape}
  {:}
  {3pt}
  {}

\theoremstyle{rem}

\newtheorem*{remark}{Remark}

\newcommand{\cinfty}{C^{\infty}}
\newcommand{\Open}{\mscr{O}}
\newcommand{\Closed}{\mscr{C}}
\newcommand{\OC}{\mscr{OC}}

\begin{document}
\maketitle

\begin{abstract}

This paper gives a construction, using heat kernels, of differential forms on the moduli space of metrised ribbon graphs, or equivalently on the moduli space of Riemann surfaces with boundary.  The construction depends on a manifold with a bundle of Frobenius algebras, satisfying various conditions.  These forms satisfy gluing conditions which mean they form an open topological conformal field theory, i.e.\ a  kind of open string theory.

If the integral of these forms converged, it would yield the purely quantum part of the partition function of a Chern-Simons type gauge theory.  Yang-Mills theory on a four manifold arises as one of these Chern-Simons type gauge theories.

\end{abstract}


\section{Introduction}

\begin{definition}
An \emph{elliptic space} $(M,\A)$ consists of 
\begin{itemize}
\item
A compact $\cinfty$ manifold $M$. Let $\cinfty(M) = \cinfty(M,\C)$ denote the algebra of complex valued smooth functions on $M$.
\item
A $\Z/2$ graded algebra $\A$ over $\cinfty(M)$, which is projective as a $\cinfty(M)$ module.     This is the same data as a finite dimensional bundle of complex $\Z/2$  graded associative algebras $\A_{M}$ on $M$, whose algebra of global sections is $\A$.
\item
An odd derivation $Q : \A \to \A$, which satisfies $Q^2 = 0$, and which is an order one differential operator with respect to $\cinfty(M)$.
\item
$(\A,Q)$ is required to be an elliptic complex.  
\end{itemize}

\end{definition}

\begin{definition}
A \emph{Calabi-Yau} structure on an elliptic space is a trace map $\Tr : \A \to \C$, which factors as
$$
\A \xto{\Tr_M} \dens M \xto{\int_M} \C
$$
where the first arrow is $\cinfty(M)$ linear. The trace map must satisfy
\begin{align*}
\op{Tr} ([a,a'] ) &= 0 \\
\op{Tr} ( Qa ) &= 0 
\end{align*}
There  is a  map $\A_M \otimes \A_M \to \op{Densities}_M$ of vector bundles on $M$, defined by $a \otimes a' \mapsto \Tr_M (aa')$.  We require that this be non-degenerate on each fibre of the vector bundle $\A_M$.
\end{definition}
The name elliptic space was suggested by Maxim Kontsevich, who studied these concepts several years ago.

The trace map $\Tr : \A \to \C$ can be even or odd, depending on whether the map $\A\to \dens{M}$ is an even or odd map of super vector bundles.  We will refer to this as the parity $p(\A)$ of $\A$.

These definitions can be modified to deal with $\Z$ graded elliptic spaces; the grading conventions are such that $Q$ is of degree $+1$.   A $\Z$ graded CY elliptic space of dimension $d$ is one such that the trace map $\Tr : \A \to \C$ is of degree $-d$.  The most interesting examples are of dimension $3$.  

   We will use the notation $\A \otimes \A$ to refer to the space of sections of the vector bundle $\A_{M} \boxtimes \A_{M}$ on $M^{2}$.  In other words, $\otimes$ will refer to the completed projective tensor product where appropriate.

\begin{definition}
Let $(M,\A)$ be a CY elliptic space.  A  Hermitian metric $\ip{\quad}$ on $\A$ is called compatible if there exists a complex antilinear, $\cinfty(M,\R)$ linear, operator $\ast : \A \to \A$ such that 
$$
\Tr (a \ast b ) = \ip{a,b}
$$
and 
$$\ast \ast \alpha = (-1)^{\abs{\alpha} (p(\A) + 1) }\alpha$$
where $\alpha \in \A$ is of parity $\abs{\alpha}$. 
We will sometimes refer to a compatible metric on $\A$ as a metric on the CY elliptic space $(M,\A)$.
\end{definition}
If $\A$ is endowed with a compatible Hermitian metric, let $Q^{\dag}$ be the Hermitian adjoint to $Q$, and let 
$$
H = [Q,Q^{\dag}] = QQ^{\dag} + Q^{\dag} Q
$$
be the Hamiltonian. $H$ is an elliptic operator, self-adjoint with respect to the Hermitian inner product, and with non-negative eigenvalues. Also the $L^{2}$ completion of $\A$ has a Hilbert space basis of eigenvectors of $H$.  As $H$ is elliptic, these eigenvectors   are in $\A$.

\begin{lemma}
\begin{enumerate}
\item
$Q^{\dag}$ and $H$ are both self adjoint with respect to the trace pairing $\Tr(ab)$ on $\A$.  
\item
There is a direct sum decomposition
$$\A = \Im  Q \oplus \Ker H \oplus \Im Q^{\dag}$$
\item
Also $\Ker H = \Ker Q \cap \Ker Q^{\dag}$.
\end{enumerate}
\end{lemma}
\begin{proof}
The first statement follows from the condition that $\ast \ast \alpha = (-1)^{\abs{\alpha} (p(\A) + 1) }\alpha$, and the second two statements follow from the ellipticity of $H$.
\end{proof}

\vspace{5pt}

\section{Summary of the paper}
  
\subsection{Examples of Calabi-Yau elliptic spaces}

There is a natural Chern-Simons type action associated to an odd CY elliptic space, given by the formula
$$
S(a)  = \Tr \left( \frac{1}{2} a Q a + \frac{1}{3} a^{3} \right) 
$$
This is a functional on the linear super-manifold $\Pi \A$. 

In section \ref{section examples}, I give some examples of Calabi-Yau elliptic spaces, and discuss the associated Chern-Simons type actions.  The most obvious example is when $M$ is a compact oriented manifold, $E$ is a vector bundle on $M$ and $\A = \Omega^\ast(\End(E))$.    This yields ordinary Chern-Simons theory.  The first non-trivial construction of this note is to give a Calabi-Yau elliptic space whose associated Chern-Simons type action is a version of Yang-Mills theory on a compact oriented four manifold with a conformal class of metrics.  

\subsection{Forms on moduli space}

In section \ref{section forms moduli}, I give a construction of forms on the moduli space of metrised ribbon graphs, associated to a Calabi-Yau elliptic space with a metric.   If $\Gamma(n,m)$ denotes the space of metrised ribbon graphs with $n$ incoming and $m$ outgoing external vertices, I define maps
$$
K : \A^{\otimes n} \to \A^{\otimes m}\otimes \Omega^{\ast}(\Gamma(n,m), \det^{\otimes p(\A)} )
$$
Here $\det$ is a certain rank one local system on $\Gamma(n,m)$.  These maps are compatible with the differentials, with gluing of ribbon graphs along external vertices, and with disjoint union.    The construction is related to Witten's construction \cite{Wit95} of a measure on moduli space associated to Chern-Simons theory.

Let $\Gamma_{g,h} \subset \Gamma(0,0)$ be the subspace of ribbon graphs of genus $g$ with $h$ boundary components.  This space is homeomorphic to the moduli space $\mc M_{g}^{h}$ of Riemann surfaces of genus $g$ with $h$  with unparameterised boundary components\footnote{Normally, the ribbon graph moduli space is viewed as  homeomorphic to the moduli space of surfaces with punctures, and a label in $\R_{> 0}$ at each puncture. However, the latter space is homeomorphic to the space of surfaces with boundary. The labels at the punctures correspond to lengths of boundary components.  A proof is given in lemma \ref{lemmargboundary}. }.   Under this homeomorphism the local system $\det$ corresponds to the orientation local system on the manifold $\mc M_{g}^{h}$.  

More generally, let $\Gamma_{g,h,n} \subset \Gamma(n,0) $ be the space of such ribbon graphs which have $n$ external vertices, with the property that the lengths of the edges attached to the external vertices is $0$, and that the external vertices are a non-zero distance apart.  Then $\Gamma_{g,h,n}$ is homeomorphic to the space $\mc M_{g}^{h,n}$ of Riemann surfaces $\Sigma$ of genus $g$, with $h$ boundary components, and with $n$ distinct ordered points in $\partial \Sigma$.

Thus, the map $K$ we construct  gives a map 
$$
K_{g,h} : \A^{\otimes n} \to \Omega^{\ast}(\mc M_{g}^{h,n}, \Or^{\otimes p(\A)})
$$
where $\Or$ refers to the orientation sheaf. 

Taking $n = 0$, we simply get a closed form $K \in \Omega^{\ast}(\mc M_{g}^{h}, \Or^{\otimes p(\A)})$.  When $p(\A) = 0$, the cohomology class of this is the same as that obtained by applying Kontsevich's construction \cite{Kon94} to the $H^\ast(\A)$ with the structure of cyclic $A_\infty$ algebra coming from the homological perturbation lemma.

In the case when $\A$ is a $\Z$ graded CY elliptic space of dimension $d$, so that the trace is of degree $-d$, then $K_{g,h}(a_1, \ldots, a_n)$ is of degree $\sum \abs{a_i} + d ( 2g-2+h)$.  Note that $\dim \mc M_{g}^{h,n} = 6g-6+3h+n$.  This makes it clear that the case $d = 3$ plays a special role.

\subsection{Chern-Simons type gauge theories}
In section \ref{section cs}, I describe a relationship between these forms and the Chern-Simons type action associated to $\A$, when $p(\A)$ is odd. The relationship is very similar to that between open string theory and Chern-Simons theory discovered by Witten \cite{Wit95}.     

Consider the algebra $\Mat_{N}(\A) = \A \otimes_{\C} \Mat_{N}$, where $\Mat_{N}$ is the algebra of $N \times N$ complex matrices.  We can think of $\Mat_{N}(\A)$ as $N \times N$ matrices with coefficients in $\A$. The operators $Q,Q^{\dag}$ extend to $\Mat_{N}(\A)$ by acting on the entries of a matrix.  Combining the trace map $\Tr : \A \to \C$ with matrix trace gives a trace map $\Tr : \Mat_{N}(\A) \to \C$.  Consider the action
$$
S(B) =  \frac{1}{2} \Tr B Q B + \frac{1}{3} \Tr B^{3}
$$
for $B \in \Mat_{N}(\A)$. Here we work with $M$ with $p(\A) = 1$, i.e.\ the trace map $\Tr : \A \to \C$ is odd.  We will consider $S$ to be an even functional on the supervector space $\Pi \Mat_{N}(\A)$.    

Let 
$$S_k(B) = \frac{1}{2} \Tr B Q B $$
be the purely quadratic part of $S$.   Let $a \in \mc H = \Ker H$, the space of harmonic elements in $\A$. We can formally write down the partition function
$$
Z_{CS}(a,\lambda,N) =  \int_{B \in \Pi \Im Q^{\dag} } e^{S(B + a \otimes \op{Id})/ \lambda} \left/   \int_{B \in \Pi \Im Q^{\dag}} e^{S_k(B)/ \lambda}   \right.
$$
$\Pi$ denotes parity change, so the integral is over the linear supermanifold $\Pi \Im Q^{\dag} \subset \Pi \Mat_{N}(\A)$.  This integral is to be understood formally, in terms of its Feynman diagram expansion.     We perform the Feynman diagram expansion complex linearly. This means that the propagator is an element of (some completion of) $\Pi \Im Q^\dag \otimes_{\C}  \Pi \Im Q^\dag$, the interaction is a $\C$-linear  map $\Pi \Im Q^\dag  \otimes_{\C} \Pi \Im Q^\dag  \otimes_{\C} \Pi \Im Q^\dag  \to \C$, etc.  

We could alternatively pick a real subspace $\mc F_{\R} \subset \Pi \Im Q^\dag$, such that $\mc F_{\R} \otimes_{\R} \C = \Pi \Im Q^\dag$,  perform the integral over the space $\mc F_\R$.  In this case we should use the $\R$-linear  Feynman diagram expansion.  Then it is easy to check we get the same answer as we did using the $\C$-linear Feynman expansion above.

The Feynman diagram expansion of $Z_{CS}(a,\lambda,N)$ is a formal expression of the form
$$
Z_{CS}(a,\lambda,N) = \exp \left( \sum_{\gamma} \lambda^{-\chi(\gamma)}  \frac{1}{\# \op{Aut}(\gamma)} w_\gamma(a) \right)
$$
where the sum is over connected trivalent ribbon graphs $\gamma$ with some external edges.   The expression $w_\gamma(a)$ denotes a certain integral over $\R_{> 0}^k$, which generally \emph{does not} converge.   Here $k$ is the number of internal edges in $\gamma$.     To make numerical sense of $Z_{CS}(a,\lambda,N)$, we need to regularise these integrals in some way.  
  
Also, we could try to integrate the form $K_{g,h}(a^{\otimes n}) \in \Omega^{\ast}(\mc M_{g}^{h,n}, \Or)$, for $a \in \mc H = \Ker H$. Top dimensional cells in $\mc M_g^{h,n}$ are labelled by trivalent ribbon graphs $\gamma$ with $n$ external vertices.  Thus, we can write
$$
\int_{\mc M_{g}^{h,n}} K_{g,h}(a^{\otimes n}) / n!  = \sum_{\gamma} \frac{1}{\op{Aut}(\gamma)} \int_{\Met_0(\gamma)} K_{g,h}(a^{\otimes n}) 
$$
where the sum is over certain trivalent ribbon graphs, and $\Met_0(\gamma)$ denotes the space of metrics on $\gamma$ which give zero length to the $n$ external edges.

It turns out that
$$
w_\gamma(a) =  \int_{\Met_0(\gamma)} K_{g,h}(a^{\otimes n}) 
$$
Both sides of this equality represent possibly non-convergent finite dimensional integrals; the integrals are of the same integrand over the same space.  

The result of this discussion is a formal equality of functions of $a \in \mc H = \Ker H$,
$$
Z_{CS}(a,\lambda,N) =   \exp \left (\sum_{\substack {g,n \ge 0, h > 0 \\ 2g-2+h + \frac{ n}{2 } > 0}} \lambda^{2g-2+h} N^{h} \int_{\mc M_{g}^{h,n}} K_{g,h}(a^{\otimes n}) / n!  \right).
$$

\subsection{Regularisation and Sen-Zwiebach formalism}
As I mentioned earlier, the integral of our forms on moduli space does not seem to converge in general.    This is because the forms have singularities on the boundary of a natural compactification of the moduli space of surfaces with boundary.  This is the compactification considered by Liu \cite{Liu02}.  In this compactification one allows surfaces with two types of singularities: nodes on the boundary or nodes in the interior.  Also one allows boundary lengths to shrink to zero. Our forms only have singularities on surfaces with nodes in the interior, and on  surfaces which have a boundary component of length zero.

One can restrict the region over which we integrate to a compact subset of the region where the forms are non-singular,  and then of course the integral converges.  One natural way to do this is to remove the locus of surfaces which have a very short\footnote{Short means in the canonical hyperbolic metric on the surface.} closed loop.      Then we find that we are integrating over part of a set of Sen-Zwiebach \cite{SenZwi96,Zwi98} open-closed string vertices.  This is what we do if we want to construct the open-closed string partition function.      

What does this construction correspond to in terms of the Chern-Simons theory?   Restricting in this way the region in moduli space over which we integrate corresponds to disallowing ``small'' subgraphs (which are not trees).  This is a very natural way of regularising Chern-Simons theory: the singularities arise precisely from such graphs.    This suggests an intriguing connection between renormalisation and the Sen-Zwiebach BV formalism.

\section{Examples}
\label{section examples}
\subsection{Chern-Simons theory}
Let $M$ be a compact oriented manifold, and let $\A = \Omega^{\ast}(M,\C)$.  The operator $Q$ is the de Rham differential $\d$, and the trace map $\Tr : \A \to \C$ is simply integration on $M$.  If we pick a metric on $M$, we get a compatible Hermitian metric on $\A$, in the usual way.  Then $Q^\dag = \d^\ast$, and $H = \tr$. 

More generally, let $E$ be a flat complex vector bundle $M$.   Then let $\A = \End(E) \otimes_{\cinfty(M)} \Omega^{\ast}_{M}$.   A combination of vector bundle trace and the usual volume element gives $M$ a volume element.  Then, as above, $Q$ is the de Rham differential on $\A$, coupled as usual to the flat connection on $\End(E)$. If we pick a Hermitian metric on $E$ and a metric on $M$, we get a metric on $\A$ which is compatible. 

The Chern-Simons type action in this example is simply the usual Chern-Simons action. When $\dim M = 3$, this has been analyzed in \cite{AxeSin92, Kon94, Wit95}.  The papers \cite{Kon94, Sch05} study also the higher dimensional generalisation.

\subsection{Holomorphic Chern-Simons}
Let $M$ be a compact complex manifold of dimension $n$, with a holomorphic volume form $\Vol_{M}$.  Let $E$ be a complex of holomorphic vector bundles on $M$.  Let $\A = \End(E) \otimes_{\cinfty(M)} \Omega_{M}^{0,\ast}$, i.e.\ $\A$ is the Dolbeaut resolution of the sheaf of holomorphic maps $E \to E$.  The natural differential on the complex $\A$ defines the operator $Q$.  The map 
$$
\A \xrightarrow{\op{Trace}_{E}} \Omega^{0,\ast}_{M} \to \Omega^{0,n}_{M} \xto{Vol_M \wedge } \Omega^{n,n}_{M}
$$
defines the trace map $\Tr : \A \to \C$.  If we pick a metric on $M$, and a Hermitian metric on each term of the complex $E$,   there is an induced compatible Hermitian metric on the $\A$.

The action in this case is the holomorphic Chern-Simons action, as studied in \cite{Wit95,DonTho98, Tho00}.

\subsection{Yang-Mills}
Let $M$ be a compact oriented $4$ manifold with a conformal structure.   Let $\Omega^{2}_{+}(M,\R)$ be the space of self dual real two forms on $M$, and let $\Omega^{2}_{+}(M) = \Omega^{2}_{+}(M,\R) \otimes_{\R}\C$.   Let $\Omega^{k}(M)$ denote the space of complex forms on $M$.     Define $\A$ to be the complex 
$$
\Omega^{0}(M) \xto{\d} \Omega^{1}(M) \xto{\d_{+}} \Omega^{2}_{+}(M)
$$
This is an elliptic space, but it is not Calabi-Yau.   Similarly, if we have a vector bundle $E$ with an anti-self dual connection $A$, the algebra 
$$
\Omega^{0}(\End(E)) \xto{\d_A} \Omega^{1}(\End(E))  \xto{\d_{A+}} \Omega^{2}_{+}(\End(E)) 
$$
defines an elliptic space.
(This construction was explained to me by Maxim Kontsevich).

\vspace{5mm}

Next we will show how to construct Calabi-Yau elliptic spaces over conformal $4$ manifolds. Let $E$ be a complex vector bundle on $M$ with a connection $A$, which satisfies the Yang-Mills equation
$$
\d_{A} F(A)_{+} = 0
$$
Here 
$$F(A)_{+} \in \Omega^{2}_{+}(\End(E)) = \Omega^2_+ (M) \otimes_{\cinfty(M)} \End(E)$$ denotes the self-dual part of the curvature. 

We can make a non-commutative elliptic space using $E$.  It is easiest to describe this construction in several steps.   First, consider the algebra $\Omega^\ast(\End(E))$, with the odd derivation $\d_A$.  Of course $\d_A^2$ is not necessarily zero.  Now take the algebra 
$$
\mc B = \Omega^\ast(\End (E)) \otimes \C[\epsilon]
$$
where $\epsilon$ is a parameter of degree $-1$.  The operator $\d_A$ extends to a derivation of $\mc B$, in a unique way such that $\d_A \epsilon = 0$.

For an element $x\in \mc B$, let $\mc L_x$ be the inner derivation given by $\mc L_x(y) = [x,y]$. Observe that for any derivation $D$ of $\mc B$, $[D,\mc L_x] = \mc L_{Dx}$.  

Define a derivation $Q$ of $\mc B$ by
$$
Q  = \d_A  -  \mc L_{\epsilon F(A)_+}   + \frac{\partial }{\partial \epsilon} 
$$
In other words, if $a,b \in \Omega^*(\End(E))$, then
$$
Q(a +  \epsilon b) = \d_A a + b -  \epsilon [F(A)_+,a]  -  \epsilon \d_A b 
$$
It is easy to calculate that 
$$[Q,Q] = 2\mc L_{F(A)_-}.$$  Indeed, the fact that $\d_A F(A)_+ = 0$ implies that $[\d_A,\mc L_{\epsilon F(A)_+} ] = 0$, the second two terms in $Q$ commute with themselves, and 
$$
[\d_A,\d_A] - 2[\frac{\partial}{\partial \epsilon}, \mc L_{\epsilon F(A)_+}] = 2\mc L_{F(A)} - 2\mc L_{F(A)_+}
$$

Next, consider the subalgebra $\mc B^+$ of $\mc B$ spanned by elements $a + \epsilon b$ such that 
$$b \in \Omega^2_+(\End(E)) \oplus \Omega^3(\End(E)) \oplus \Omega^4(\End(E))$$
The operator $Q$ preserves the subalgebra $\mc B^+$.

Now let $I \subset \mc B^+$ be the ideal
$$ I = \Omega^2_-(\End(E)) \oplus \Omega^3(\End(E)) \oplus \Omega^4(\End(E)) \subset \mc B^+$$
The ideal $I$ is preserved by $Q$.  Therefore we can define the quotient
$$
\A = \mc B^+ / I
$$
In the original algebra $\mc B$, $[Q,Q]=\mc L_{F(A)_-}$. It follows that $[Q,Q] = 0$ in $\A$.

The differential graded algebra $\mc A$ looks like 
$$
\xymatrix{   \Omega^{0} (\End(E)) \ar[dd]_{\d_{A}} \ar[ddrr]^{-[F(A)_{+},\quad]} &  &  &  \text{ degree } 0  \\ 
 & & & \\
 \Omega^{1}(\End(E)) \ar[dd]_{\d_{A+}}   \ar'[dr]|{-[F(A)_{+}, \quad]} [ddrr] & \oplus &  \Omega^{2}_{+} (\End(E)) \ar[dd]^{- \d_{A}} \ar[ddll]|(.3){\op{Id}} &  \text{ degree } 1    \\
 & & & \\
  \Omega^{2}_{+}(\End(E))  \ar[ddrr]^{-[F(A)_{+},\quad]} & \oplus & \Omega^{3} (\End(E)) \ar[dd]^{-\d_{A}} &  \text{ degree } 2  \\
  & & & \\
   &  & \Omega^{4} (\End(E)) &  \text{ degree } 3  }
$$
The signs are slightly tricky in this diagram. We should think of the right hand column as being things of the form $ \epsilon b$, where $b \in \Omega^\ast (\End(E))$.  Therefore the left action of the left column on the right column is twisted by a sign.  The right action is the obvious one. This explains also why we have $-\d_A$ instead of $\d_A$ in the right hand column.

The trace map $\Tr : \A \to \C$ is given by the composition of the  trace map $\Tr_E : \Omega^{4} (\End(E)) \to \Omega^{4}(M)$ with integration on $M$. 

\begin{lemma} 
$(M,\A)$ is a Calabi-Yau elliptic space.   Also, if we pick a metric on $M$ in the given conformal class, and a Hermitian metric on $E$, we get a compatible Hermitian metric on the space $\A$, and so a metric on the elliptic space $(M,\A)$.
\end{lemma}
\begin{proof}
We have already seen that $\A$ is associative and  $Q$ is a derivation of square zero. It is also easy to see that $\A$ is elliptic, and that the pairing $\Tr(ab)$ is non-degenerate.    

 Let $\Tr_{E} : \Omega^{*}(\End(E)) \to \Omega^{\ast}_{M}$ denote the fibrewise trace map. For any connection $A$ on $E$, and any $x \in \Omega^{3}(\End(E))$, $\Tr (d_{A} x) = \d \Tr (x)$.   Also, for $x,x' \in \Omega^2(\End(E))$, we have $\Tr_E([x,x']) = 0$.  It follows that for $a \in \A^3$, $\Tr Qa = 0$.   
 
Picking a metric on $M$ and $E$ gives complex antilinear Hodge star maps
 $$
 \ast_0 : \Omega^i(\End(E)) \to \Omega^{4-i}(\End(E))
 $$
 defined by 
 $$
 \ast_0 ( \omega \otimes X) =  \br{\ast(\omega)} \otimes X^\dagger
 $$
 where $X \in \End(E)$ and $\omega \in \Omega^i(M)$.  $X^\dagger$ denotes the Hermitian adjoint to $X$ and $\ast(\omega)$ is the usual Hodge star of $\omega$.
 
The Hodge star operator $\ast$ on $\A$ as follows.  As always, we are thinking of $\A$ as a subquotient of $\Omega^\ast(\End(E)) [\epsilon]$. Then $\ast$ is defined by
$$
\ast (  a+ \epsilon b ) = (-1)^{\abs{a}} \epsilon  \ast_0 a +  \ast_0 b
$$
 It is  clear that $\ast^2 =  1$.  The pairing $\Tr (a \ast a')$ is the Hermitian metric on each direct summand of $\A$ which comes from the choice of  metric on $M$ and on $E$.  This is obvious for each summand except $\Omega^1(\End E) \subset \A^1$. But, if $a,a'$ are in this space, then 
$$\Tr ( a \ast a' ) = - \Tr ( a \ast_0 a' )  = \Tr ( (\ast a' ) a )$$ 
as $\Tr ( [x,y]) = 0$.  

\end{proof}

\subsubsection{The action functional}

The Chern-Simons type action functional, for $\alpha \in \A^{1}$, is defined by
$$
S(\alpha) = \frac{1}{2}\Tr \alpha Q \alpha + \frac{1}{3} \Tr \alpha^{3}
$$
Let us write $\alpha = a + \epsilon b$, where $a \in \Omega^1(\End(E))$, $b \in \Omega^2_+(\End(E))$, and $\epsilon$  is an odd parameter.  As before, we are considering $\A$ as a subquotient of the algebra $\Omega^\ast(\End(E)) [\epsilon]$.
\begin{proposition}
$$
S ( a + \epsilon b ) =  \Tr  \left( \frac{1}{2} b^2  + (F(A + a) - F(A)) (b - F(A)_+) \right)
$$
\end{proposition}
\begin{proof}
\begin{align*}
Q \alpha  &=  \d_{A+} a + b  - \epsilon \d_A b - \epsilon [F(A)_+, a]\\
\alpha Q \alpha &=  \epsilon b \d_A a + \epsilon b^2 +  \epsilon a \d_A b  + \epsilon a [F(A)_+, a]  \\
 \alpha^3 &=  \epsilon a^2 b  + \epsilon b a^2 - \epsilon a b a 
\end{align*}
so that 
\begin{align*}
S(a + \epsilon b) &= 
\Tr \left ( \frac{1}{2} b^2  +  (\d_A a) b + a^2 b + a [F(A)_+, a] \right) \\
&= \Tr  \left( \frac{1}{2} b^2  + (F(A + a)_+ - F(A)_+ ) b - a ^2 F(A)_+ \right) 
\end{align*}
 where $A + a$ refers to the connection obtained by adding $a$ onto the given Yang-Mills connection $A$ on $E$.   Since 
$$(\d_A a + a^2 ) F(A)_+ = (F(A+a)_+ - F(A)_+ ) F(A)_+$$
and 
$$
\Tr (( \d_A a) F(A)_+ ) = 0
$$
because $\d_A F(A)_+ = 0$, this action can be rewritten
\begin{align*}
S(a + \epsilon b ) &= \Tr  \left( \frac{1}{2} b^2  + (F(A + a)_+ - F(A)_+ ) b - F(A + a)_+ F(A)_+  + F(A)_+^2  \right)  \\
&=  \Tr  \left( \frac{1}{2} b^2  + (F(A + a)_+ - F(A)_+ ) (b - F(A)_+) \right)
\end{align*}
as desired. 
\end{proof}
\begin{lemma}
The set of critical points of the action $S$ is the set $a + \epsilon b$ with $\d_{A + a} F(A+a)_+ = 0$ and $b = F(A)_+ - F(A+a)_+$.
\end{lemma}
\begin{proof}
Varying $b$ immediately leads to the equation that at a critical point, 
$$
b =     F(A)_{+} -  F(A+a)_{+}
$$
If we vary $a$ to $a + \delta a$, then, to first order,
$$
F(A + a + \delta a ) = F(A + a) + \d_{A+a} \delta a
$$

Therefore we find the constraint
$$
\d_{A+a} (b -  F(A)_{+}) = 0
$$
Finally, since at a critical point $b = F(A)_+ - F(A+a)_{+} $, this equation becomes the Yang-Mills equation
$$
\d_{A+a} F(A+a)_{+} = 0
$$
\end{proof}

Thus we have shown that the set of critical points is the same as the set of those $a$ satisfying the Yang-Mills equation.  In other words, the action we wrote down is classically equivalent to the Yang-Mills action.

\subsubsection{Relationship to Yang-Mills at the quantum level}
At the quantum level, the Chern-Simons type theory we construct is equivalent to ordinary Yang-Mills, given by the action
$$
S_{\YM}(a) =  \frac{1}{2} \Tr ( F(A + a)_+^2 )
$$
This equivalence seems to be well-known to experts in the area.   

If we perform the simple change of coordinates $b \mapsto b + F(A)_+$, then our action becomes
$$
\Tr  \left( \frac{1}{2} b^2   + F(A + a)_+ b + \frac{1}{2} F(A)_+^2   \right)
$$
Since $F(A)_+^2$ is a constant, we can ignore it.  Then the action is precisely the same as that studied in \cite{Wit03}, section 4.4, where it was described as equivalent to ordinary Yang-Mills.

A very similar action was studied in \cite{CatCotFuc98,MarZen96}, given by the formula
$$\Tr  \left( \frac{1}{2} b \ast b   + F(A + a) b   \right).$$
Here we are perturbing around a flat connection $A$,  $a \in \Omega^1(\End(E))$, and $b \in \Omega^2(\End(E))$ is not necessarily self dual. These authors showed this theory was equivalent to Yang-Mills at the quantum level. If we perturb around a connection $A$ with $F(A)_+ = 0$, then our action is given by the same formula. The only difference is that our field $b$ is self dual.

Besides these references to the literature, let me give a simple formal argument showing that the theory we considered here should be  equivalent to ordinary Yang-Mills.   Define a non-linear map $\Phi : \A^1 \to \A^1$ by
$$
\Phi ( a + \epsilon b  ) =a + \epsilon ( b + F(A)_+ -  F(A+a)_{+}) $$
It is easy to see that 
\begin{align*}
S( \Phi (a + \epsilon b) ) &= \frac{1}{2} \Tr ( F(A + a)_+^2 - F(A)_+^2  + b^2)   \\
&= - S_{\YM} (a ) + \frac{1}{2} \Tr b^2 - \frac{1}{2} \Tr F(A)_+^2
\end{align*}
The last term in this formula is a constant, and can be ignored. The derivative of the map $\Phi$ is the identity plus a triangular matrix, and therefore the Jacobian is trivial.   

If we consider the theory defined by the action $- S_{\YM} (a ) + \frac{1}{2} \Tr b^2$, and we integrate out $b$, we recover the usual Yang-Mills action.  

Let $\alpha \in \Omega^3(\End(E))$, considered as a linear form on $\Omega^1(\End(E))$. The argument above show that, formally,  we have the identity
\begin{align*}
Z(\alpha) &= \int_{a,b} \exp \lambda^{-1} \left(  \Tr  \left( \frac{1}{2} b^2  + (F(A + a) - F(A)) (b - F(A)_+) \right)  + \alpha (a) \right)  \\
&= c  \int_a \exp  \lambda^{-1} \left(  -S_YM(a) + \alpha(a) \right)   \int_b \exp \lambda^{-1} \left( \frac{1}{2} \Tr b^2  \right)\\
&= c' \int_a \exp  \lambda^{-1} \left(  -S_YM(a) + \alpha(a) \right)
\end{align*}
where $c,c'$ are constants independent of $\lambda$. 

\subsubsection{Categorical generalisation}
This construction can be easily generalised to give a differential graded Calabi-Yau category associated to the conformal $4$ manifold $M$,  whose objects are vector bundles with Yang-Mills connections.    

If $E$ is such a vector bundle, let us use the notation $\Hom_{\YM}(E,E)$ for the differential graded algebra $\A$ constructed above.  If $E_1,E_2$ are Yang-Mills bundles, the algebra $\End(E_1 \oplus E_2)$ has two idempotents $\pi_1,\pi_2$, corresponding to projection onto $E_1$ and $E_2$.  These extend to idempotents in $\Hom_{\YM}(E_1 \oplus E_2, E_1 \oplus E_2)$, which commute with the differential $Q$.  Let 
$$
\Hom_{\YM}(E_1,E_2) = \pi_2 \Hom_{\YM}(E_1 \oplus E_2, E_1 \oplus E_2)  \pi_1 
$$
This is a complex which looks like
$$
\xymatrix{   \Omega^{0} (\Hom(E_1,E_2)) \ar[dd]_{\d_{A}} \ar[ddrr]^{F(A_1)_{+} -  F(A_2)_+} &  &  &  \text{ degree } 0  \\ 
 & & & \\
 \Omega^{1}(\Hom(E_1,E_2)) \ar[dd]_{\d_{A+}}   \ar'[dr]|{F(A_1)_{+} -  F(A_2)_+} [ddrr] & \oplus &  \Omega^{2}_{+} (\Hom(E_1,E_2)) \ar[dd]^{- \d_{A}} \ar[ddll]|(.3){\op{Id}} &  \text{ degree } 1    \\
 & & & \\
  \Omega^{2}_{+}(\Hom(E_1,E_2))  \ar[ddrr]^{F(A_1)_{+} -  F(A_2)_+} & \oplus & \Omega^{3} (\Hom(E_1,E_2)) \ar[dd]^{-\d_{A}} &  \text{ degree } 2  \\
  & & & \\
   &  & \Omega^{4} (\Hom(E_1,E_2)) &  \text{ degree } 3  }
$$
The notation $F(A_1)_+ - F(A_2)_+$ denotes the operation $x \mapsto x F(A_1)_+   - F(A_2)_+ x$.  

This defines the morphism space in the differential graded Calabi-Yau category of Yang-Mills bundles on $M$.   Passing to the cohomology category yields a Calabi-Yau $A_\infty$ category.

One can associate  to each such Calabi-Yau $A_\infty$ category  a closed TCFT, i.e. a kind of closed string theory (\cite{Kon03, Cos04, KonSoi06}).  The closed string states is the Hochschild homology of the category.   It seems to me to be a very interesting problem (which I have no idea how to solve)  to understand in geometric terms the associated closed TCFT.

\subsection{$G2$ manifolds}
There is also a construction of Calabi-Yau elliptic spaces associated to $G2$ manifolds.  This construction was kindly explained to me by Maxim Kontsevich; it is also discussed in the article \cite{DonTho98}.   

Let $M$ be a compact $G2$ manifold.  Let $\sigma \in \Omega^4(M)$ be the canonical $4$ form.  Let $E$ be a complex vector bundle on $M$ with a connection $A$ such that $\sigma \wedge F(A) = 0$.  Let $\op{Ann}(\sigma \otimes \op{Id}_E) \subset \Omega^\ast(\End(E))$ be the ideal of elements whose  product with $\sigma \otimes \op{Id}_E$ is zero. Let 
$$
\A = \Omega^\ast(\End(E)) / \op{Ann}(\sigma \otimes \op{Id}_E)
$$
$\A$ is concentrated in degrees $0,1,2,3$.   The covariant differentiation operator $\d_A$ descends to $\A$, because $\d_A (\sigma \otimes \op{Id}_E) = 0$. The operator $Q$ is defined to be $\d_A$ on $\A$.  The condition that $\sigma \wedge F(A) = 0$ implies that $Q^2 = 0$.

The trace map $\Tr : \A \to \C$ comes from 
$$
\Omega^3(\End(E)) \xto{\Tr_E} \Omega^3(M) \xto{\sigma \wedge} \Omega^7(M) \xto{\int_M} \C
$$
The action $S$, on an element $a \in \A^1 = \Omega^1(\End(E))$, is given by the formula
$$
S(a) = \int_M   \sigma \wedge \frac{1}{2} \Tr_E a \d_A a + \sigma \wedge \frac{1}{3} \Tr_E  a^3
$$

\subsection{$GL(n,\C)$ versus $U(n)$}
In each of the examples of Calabi-Yau elliptic spaces above, we have an associated gauge theory, given by the Chern-Simons type action.  One can ask what the gauge group is.  Let's discuss this in the simplest case, that of a flat bundle $E$ on a compact $3$ manifold $M$, with $\A = \Omega^\ast(\End(E))$.  The bundle $E$ is a flat $GL(n,\C)$ bundle.  However, when we compute the perturbative Feynman integral, we perform the Feynman rules complex linearly.  This is equivalent to integrating over \emph{any} real slice of $\Pi \Im Q^\dag$, that is any real subspace $\Pi \Im Q^\dag_\R$ such that $\Pi \Im Q^\dag_\R \otimes_\R \C = \Pi \Im Q^\dag$.    If $E$ is a flat $U(N)$ bundle, so that $E$ has a Hermitian metric compatible with the connection,  then we can look at the subspace 
$$\A_\R  = \Omega^\ast(M,\R) \otimes \text{Hermitian endomorphisms of } E \subset \A$$
Then $\Pi \A_\R \cap \Pi \Im Q^\dag$ is a real slice of $\Pi \Im Q^\dag$, and we find we are computing perturbative $U(n)$ gauge theory.

A similar discussion holds in the Yang-Mills and $G2$ cases, when the connection on the vector bundle $E$ is compatible with a Hermitian metric.    In the holomorphic Chern-Simons case, I don't see how to pick out a natural real slice.

\section{Constructing forms on the moduli space of ribbon graphs}
\label{section forms moduli}

The  construction comes from a kind of supersymmetric quantum mechanics on the non-commutative manifold $(M,\A)$.    As is usual in quantum mechanics, points in the space are replaced by functions on it.  So the ``Hilbert space'' of the theory will be $\A$.   Functions $f \in \A$ will evolve by the Hamiltonian $H$.  Also these functions will interact; any number of functions can interact at one time.  The interaction is simply given by the product of functions. The theory will give us forms on the space of metrised ribbon graphs, or equivalently on the moduli space of Riemann surfaces.    More generally, if we consider ribbon graphs with $n$ incoming and $m$ outgoing external vertices, we find a map from $\A^{\otimes n}$ to forms on moduli space tensored with $\A^{\otimes m}$.  (As always, tensor product means the completed projective tensor product).

The simplest part of construction yields, for each rooted, metrised ribbon tree $T$ with $n$ leaves, a map $\A^{\otimes n} \to \A$.  Before I describe the construction in general, I will give an informal description of this part of the theory.

\begin{figure}

\includegraphics[ bb = 1 -10 200 200]{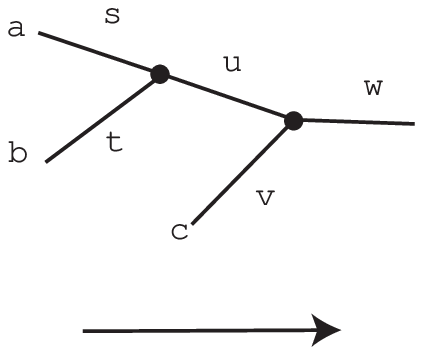}
\caption{
In this diagram, $s,t,u,v,w \in \R_{> 0}$ and $a,b,c \in \A$.  The output of this diagram is 
$$
e^{-w \tr} \left(  \left(   e^{-u\tr} \left(  \left( e^{-s \tr } a \right) \left(  e^{-t \tr  } b  \right)    \right)  \right) \left( e^{-v \tr} c    \right) \right)
$$
\label{figure_tree}
}

\end{figure}

We will consider the tree as defining a diagram for the interaction of $n$ ``particles'' to yield a single particle.    At each leaf of $T$ put a function in $\A$.  If the edge connecting to this leaf is of length $t$, apply $e^{-tH}$ to the function.  If a vertex has $k$ incoming edges, each edge labelled by a function in $\A$, multiply all these functions together to yield the function associated to the start of the outgoing edge.  Then, to get the function associated to the end of this outgoing edge, apply $e^{-l(e) H}$.  Repeating this procedure for all vertices and edges associates to our $n$ incoming functions a single function in $\A$.  An example of this procedure is given in figure \ref{figure_tree}.

In the form version of the construction, the length of an edge does not lie in $\R_{> 0}$ but in the supermanifold $\R_{> 0} \times \R^{0,1}$.  We will think of this as the odd tangent bundle of $\R_{> 0}$, so that functions on this supermanifold will be identified with forms on $\R_{> 0}$.   In this situation, a function is propagated along an edge by 
$$
f \mapsto e^{- l(e) H} f - \d l(e)  Q^{\dag} e^{-l(e) H}  f \in \A \otimes \Omega^{\ast}(\R_{> 0})
$$
At the vertices we simply multiply the incoming functions, as before.

When we use general graphs instead of trees, the edges are no longer oriented.  In this situation, the analog of the construction above is to put at each edge $e$ of length $l(e)$, connecting two vertices $v_{1},v_{2}$, the kernel 
$$K_{l(e)} (x_{1},x_{2}) - \d l(e) Q^{\dag}_{1} K_{l(e)}(x_{1},x_{2})$$
where $K_{t}$ is the heat kernel for $H$; and at each external incoming vertex an input function in $\A$.  Then, at each vertex, we multiply the elements of $\A$ corresponding to each germ of edge at that vertex, and apply the trace map $\Tr : \A \to \C$.  This is done using the cyclic order at the vertex.  This yields an element of the tensor product of the space of forms on graphs, with the space $\A^{\otimes \text{outgoing vertices}}$.

\subsection{Heat kernels}

Let $K \in \A^{\otimes 2}$.   $K$ defines a convolution operator $\A \to \A$, by
$$
f \mapsto K f \defeq  (-1)^{p(\A) \abs{f}}  \Tr_{y} K(x,y)  f(y)
$$
The sign occurs so that  kernels closed under $Q_{x} + Q_{y}$ correspond to operators which commute with $Q$.    We will generally use the same symbol for a kernel and the associated operator.
\begin{definition}
A \emph{heat kernel} for $H$ is a an element 
$$
K \in \cinfty(\R_{> 0}) \otimes \A \otimes \A  = \Gamma(\R_{> 0} \times M \times M, \A_{M} \boxtimes \A_{M})
$$
such that the convolution operator $K_{t} : \cinfty(E) \to \cinfty(E)$ behaves like $e^{-t H}$, where $t$ is the coordinate on $\R_{> 0}$.  That is, 
\begin{align*}
\frac{\d}{\d t} K_{t} f &= -H  K_{t} f \\
\lim_{t \to 0}K_{t}  f &= f
\end{align*}
\end{definition}
The heat equation can be expressed more directly as 
$$
\frac{\d}{\d t} K_{t}(x,y) = -H_{x}K_{t}(x,y)
$$
and the second condition says that $\lim_{t \to 0} K_{t}(x,y) = \delta_{x,y}$ is a delta distribution on the diagonal.  

The results of Gilkey \cite{Gil95} imply that $H$ admits a heat kernel.  This is a general property  of Laplacian operators which arise from elliptic complexes.   

\subsubsection{Identities satisfied by the heat kernel}

Let us consider the heat kernel $K_{t}$ as an operator on $\A$.  Viewing the operator $K_{t}$ as $e^{-t H}$, it is obvious that $K_{t}$ is self-adjoint, and that 
\begin{align*}
[Q,K_{t}] = [Q^{\dag},K_{t}] &= 0 \\
[H,K_{t}] &= 0
\end{align*}
Let 
$$
L_{t}(x,y) = - Q^{\dag}_{x}K_{t}(x,y)
$$
The kernel $L_{t}$ represents the operator $f \mapsto -Q^{\dag} K_{t} f$.  Therefore $L_{t}$ is self adjoint and
$$[Q,L_{t}] = - H K_{t}$$ 
In kernel terms, these identities become the following.
\begin{proposition}
\label{prop heat kernel identities}
The identities
\begin{align*}
Q_{x} K_{t}(x,y) +  Q_{y} K_{t}(x,y) &= 0 \\
Q_{x}^{\dag} K_{t}(x,y) &=  Q_{y}^{\dag} K_{t}(x,y)  \\
H_{x}K_{t}(x,y) &= H_{y} K_{t}(x,y)  \\
K_{t}(x,y) &= (-1)^{p(\A)}K_{t}(y,x) \\
L_{t}(x,y) &= (-1)^{p(\A)} L_{t}(y,x) \\
(Q_{x}+Q_{y})L_{t} (x,y) &= -H_{x} K_{t}(x,y) = \frac{\d}{\d t} K_{t}(x,y)
\end{align*}
are satisfied.
\end{proposition}
The proof, which  is given in an appendix, is  easy except for some awkwardness with signs.

\subsection{Categories of metrised graphs}

\begin{figure}
\includegraphics[ bb = 1 80 200 170]{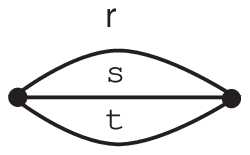}
\caption{  The theta graph. $r,s,t$ refer to coordinates on the moduli space of metrics on this graph. \label{figure_theta_graph}}
\end{figure}

Before I give the general construction of forms on the space of metrised ribbon graphs, let me give a simple example.    The theta graph is illustrated in figure \ref{figure_theta_graph}.   This space of metrics is a dense open subset of $\R_{\ge 0}^{3}$. Coordinates on the moduli space of metrics on this graph are given $r,s,t$.  Let us write, formally,
$$
K_{t} + \d t L_{t} = \sum \alpha_{i} \otimes \alpha'_{i} \otimes \beta_{i}(t) \in \A \otimes \A \otimes \Omega^{\ast}(\R_{> 0})
$$
Then the form is defined by 
$$
\sum_{i,j,k} \Tr \left(  \alpha_{i} \alpha_{j}\alpha_{k}  \right) \Tr \left(  \alpha'_{i} \alpha'_{k}\alpha'_{j}  \right)  \beta_{i}(r) \beta_{j}(s) \beta_{k}(t)
$$
(with appropriate signs).

\begin{definition}
A \emph{graph} $\gamma$ is a compact one dimensional cell complex.  The $0$ cells are called vertices, $V(\gamma)$, and the one cells the edges, $E(\gamma)$.

A \emph{ribbon graph} is a graph with a cyclic order on the set of germs of edges emanating from each vertex.
\end{definition}

Let $\Gamma(n,m)$ be the space of metrised ribbon graphs with $n$ incoming and $m$ outgoing labelled external vertices.  These external vertices are uni-valent; all other vertices are at least tri-valent. These graphs need not be connected.  

An edge is called internal if neither of its ends are an external vertex, otherwise the edge is external.   To each edge of a graph $\gamma \in \Gamma(n,m)$ is assigned a length $l(e) \in \R_{\ge 0}$.     Let $[n],[m]\subset V(\gamma)$ denote the sets of incoming and outgoing vertices. 

\begin{figure}
\includegraphics[ bb = 1 -10 200 200]{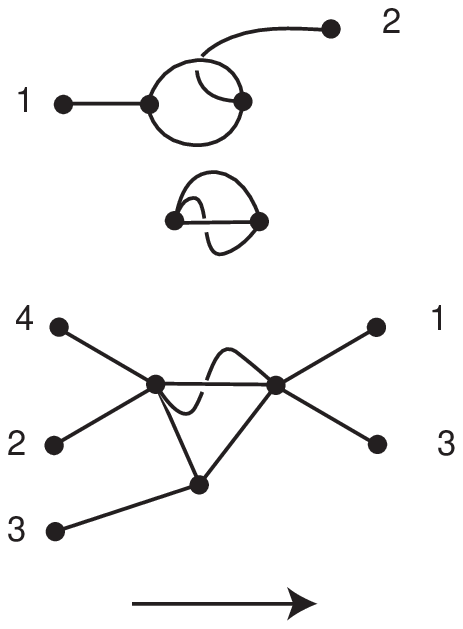}
\caption{A typical ribbon graph, with $4$ incoming vertices (on the left) and $3$ outgoing vertices.}
\end{figure}

The following two conditions are imposed on these metrised graphs.
\begin{enumerate}
\item
Every closed loop of edges in $\gamma$ is of positive length.
\item
Every path of edges in $\gamma$ which starts and ends at different outgoing external vertices is of positive length. 
\end{enumerate}
 If $e$ is an internal edge of $\gamma$ of length zero, then we identify $\gamma$ with $\gamma/e$, the graph obtained by contracting $e$. 
 
 We allow connected components of graphs to be of the following exceptional types.
 \begin{enumerate}
\item
A graph with two external vertices, one edge, and no internal vertices.  The external vertices can be any configuration of incoming or outgoing.
\item
A graph with two internal vertices, and two edges between them.  As these internal vertices are not trivalent, this violates the conditions above. 
\end{enumerate}

The space $\Gamma(m,n)$ is closely related to moduli spaces of Riemann surfaces with open boundary, as I will discuss later.

There are continuous gluing maps
$$
\Gamma(m,l) \times \Gamma(n,m) \to \Gamma(n,l)
$$
The conditions on the allowed lengths of tails are imposed so that when we glue two graphs, we never end up with a cycle of zero length.

There are disjoint union maps
$$
\Gamma(s,r) \times \Gamma(n,m) \to \Gamma(s+n,m+r)
$$
This structure makes the spaces $\Gamma(n,m)$ into the morphisms in a topological symmetric monoidal category whose objects are the non-negative integers $\Z_{\ge 0}$. The units in this category are given by graphs each of whose connected components has one incoming and one outgoing vertex, and no other vertices.   Let us denote this category by $\Gamma$.

\subsection{A local system on $\Gamma(n,m)$}

We want to take chains on $\Gamma(n,m)$ with coefficients in a certain local system.  

Let $\det$ be the graded local system on $\Gamma(n,m)$ whose fibre at a metrised graph $\gamma$ is
$$
\det(\gamma) = \Pi ^{m-\chi(\gamma) }\det(H^{\ast}(\gamma,[m])) 
$$
that is, $\det(\gamma)$ is the determinant of the relative cohomology group.  
The notation $\Pi^{m-\chi(\gamma)}$ refers to a change of parity.

We can describe this more explicitly.  A trivialisation of the line $\det$ is given by an ordering on the set  $V(\gamma)\setminus [m]$, an ordering on the set of edges $E(\gamma)$, and an orientation on each edge.  If we change the orientation of an edge we change this element by a sign.  If we change the orderings we change by a sign given by the signature of the permutation.

Under the composition map
$$
\Gamma(m,l) \times \Gamma(n,m) \to \Gamma(n,l)
$$
the pullback of the  local system $\det$ is canonically isomorphic to $\det \boxtimes \det$, in  an ``associative'' fashion.  That is, under the map for the composition of three morphisms in $\Gamma$, the two  isomorphisms between the pullback of $\det$ and $\det^{\boxtimes 3}$ coincide.  

The local system $\det$ is similarly well behaved under disjoint union maps.

Recall that there is a notion of boundary cycle in a ribbon graph, corresponding to boundaries of Riemann surfaces.    Let $\Gamma_{g,h} \subset \Gamma(0,0)$ be the subspace of graphs of genus $g$ with $h$ boundary cycles.  More generally, let $\Gamma_{g,h,n} \subset \Gamma(n,0)$ be the space of such ribbon graphs which have $n$ external vertices, with the property that the lengths of the edges attached to the external vertices is $0$, and that the external vertices are a non-zero distance apart.   We only consider these spaces when $2g-2+h + \frac{1}{2}n > 0$.

\begin{lemma}
$\Gamma_{g,h}$ is homeomorphic to the space $\mc M_g^h$ of Riemann surfaces of genus $g$ with $h$ unordered boundaries. Also $\Gamma_{g,h,n}$ is homeomorphic to the space $\mc M_{g}^{h,n}$ of Riemann surfaces $\Sigma$ of genus $g$, with $h$ boundary components, and with $n$ distinct ordered points in $\partial \Sigma$.

\label{lemmargboundary}
\end{lemma}

\begin{proof}
Let $\til {\mc M}_{g}^{h}$ be the moduli space of Riemann surfaces of genus $g$ with $h$ unordered punctures, with each puncture labelled by an element in $\R_{> 0}$.   The standard ribbon graph decomposition says that $\Gamma_{g,h}(0,0)$ is homeomorphic to $\til {\mc M}_{g}^{h}$.  It suffices to show that $\til {\mc M}_{g}^{h}$ is homeomorphic to $\mc M_{g}^{h}$. 

Each $\Sigma \in \mc M_{g}^{h}$ has a canonical hyperbolic metric with geodesic boundary. This induces a  parameterisation up to rotation on the boundary of $\Sigma$, given by arc-length.  If $\Sigma \in \mc M_{g}^{h}$, let $\Sigma'$ be the space obtained from gluing a copy of the punctured disc $\{ z \in \C \mid 0 < \abs{z} \le 1 \}$ onto each boundary component. This gluing is performed using the standard parameterisation on the boundary of the disc, and the arc-length parameterisation on the boundary of $\Sigma$.  This is well-defined as the parameterisation on the boundary of the disc is rotation invariant.    $\Sigma'$ is a surface with $h$ punctures; label these punctures by the lengths of the corresponding boundary components of $\Sigma$.

Sending $\Sigma$ to $\Sigma'$ with these labels defines a map of orbifolds $\mc M_{g}^{h} \to \til {\mc M}_{g}^{h}$.  It is easy to check that this is a homeomorphism. Indeed,   using Fenchel-Nielsen coordinates on Teichmuller space, it suffices to check the case  when $g = 0$ and $h = 3$.

The homeomorphism $\Gamma_{g,h,n} \iso \mc M_{g}^{h,n}$ can be proved by a similar argument. 

\end{proof}

\begin{lemma}
On the manifold $\Gamma_{g,h,n}$, $\det$ is naturally isomorphic to the orientation sheaf, with parity $2g-2+h$. \label{det_or}
\end{lemma}
\begin{proof}
The fibres of the map $\Gamma_{g,h,n} \to \Gamma_{g,h}$ given by forgetting all the external tails, are naturally oriented. Also the sheaf $\det$ on $\Gamma_{g,h,n}$ is pulled back from $\Gamma_{g,h}$. Therefore it suffices to prove the result for $\Gamma_{g,h}$.

We've seen above that the space $\Gamma_{g,h}$ is homeomorphic to the space of Riemann surfaces of genus $g$ with $h$ boundary components, or equivalently to the space of surfaces of genus $g$, with $h$ punctures, and a map from the set of punctures to $\R_{> 0}^{h}$.  The punctures are unordered.   

The space of surfaces with punctures is naturally oriented, as it is a complex manifold.  Therefore the orientation sheaf is naturally isomorphic to $\det (\R^{h})$.

The sheaf $\det$ can be viewed as the determinant of the cohomology of the fibre of the universal surface on the space of surfaces with punctures.  Let $\Sigma$ denote the universal surface, and $\br{\Sigma}$ be the compactification of $\Sigma$ by filling in the punctures.

The cohomology $H^{1}(\br {\Sigma})$ is naturally oriented, as it is a symplectic vector space.  It follows that the determinant line $\det H^{\ast} (\br{\Sigma})$ is naturally trivial.  

From the  sequence relating the cohomology of $\Sigma$ with that of $\br \Sigma$ we see that  $\det H^{\ast}(\Sigma)$ is naturally isomorphic to $\det \R^{h} \otimes \det H^{\ast}( (S^{1})^{h})$, which is naturally isomorphic to $\det \R^{h}$.
\end{proof}

\subsection{Cellular differential forms}

The space $\Gamma(n,m)$ has a natural decomposition into non-compact orbi-cells.  For each graph $\gamma$, without a metric, with $n$ incoming and $m$ outgoing external edges, let $\Met(\gamma) \subset \R_{\ge 0}^{E(\gamma)}$ be the space of allowed metrics on $\gamma$.    The orbi-cells of $\Gamma(n,m)$ are by definition the images of the spaces $\Met(\gamma)$.

If $e$ is an internal edge of $\gamma$ which isn't a loop, there is an embedding 
$$
\Met(\gamma/e)\into \Met(\gamma)
$$
as the subspace of metrics which give $e$ length zero.
\begin{definition}
The space $\Omega^{i}_{cell}(\Gamma(n,m))$ of cellular differential forms is defined as follows. An element $\omega\in \Omega^{i}_{cell}(\Gamma(n,m))$ consists of a form  form $\omega_{\gamma} \in \Omega^{i}(\Met(\gamma),\C)$, for each graph $\gamma$ as above, which is $\Aut(\gamma)$ equivariant, and such that, for each edge $e$ of $\gamma$ which is not a loop,
$$
\omega_{\gamma/e} = \omega_{\gamma} \mid_{\Met(\gamma/e)}
$$
\end{definition}
The complex $\Omega^{i}_{cell}(\Gamma(n,m))$, with the de Rham differential, computes the cohomology of $\Gamma(n,m)$.  

The composition maps $\Gamma(m,l) \times \Gamma(n,m) \to \Gamma(n,l)$ induce coproduct maps
$$
\Omega^{\ast}_{cell}(\Gamma(n,l)) \to \Omega^{\ast}_{cell}(\Gamma(m,l) \times \Gamma(n,m)) = \Omega^{\ast}_{cell}(\Gamma(m,l) )\otimes \Omega^{\ast}_{cell}(\Gamma(n,m) )
$$
where as usual  $\otimes$ refers to the completed projective tensor product.

Denote by $C_{*}(\Gamma(n,m))$ the subcomplex of normalised singular simplicial chains on $\Gamma(n,m)$, spanned by simplices which lie entirely in one of the closed cells, and which are $\cinfty$.  

The complexes $C_{\ast}(\Gamma(n,m))$ form the morphisms in a  differential graded symmetric monoidal category $C_{\ast}(\Gamma)$, whose objects are the non-negative integers.    

Because of the compatibility between the local system $\det$ and the composition and disjoint union maps, we can form a twisted dg category $C_{\ast}(\Gamma,\det)$ whose morphism complexes are singular chains with local coefficients $C_{\ast}(\Gamma(n,m),\det)$.

Note there is an integration pairing
$$
C_{i}(\Gamma(n,m)) \otimes \Omega^{i}_{cell}(\Gamma(n,m)) \to \C
$$
The composition on $C_{\ast}(\Gamma)$ and the coproduct on $\Omega^{\ast}(\Gamma)$ are of course dual under this pairing.  

\subsection{Forms on moduli of graphs}
Let $\gamma$ be a graph with $n$ incoming and $m$ outgoing external edges, but without a metric.   As before, let $\Met(\gamma)$ be the space of allowed metrics on $\gamma$.

Recall that $[m] \subset V(\gamma)$ is the set of outgoing external vertices.
Let us pick an ordering on the sets $V(\gamma) \setminus [m]$ and $E(\gamma)$, and an orientation of each edge.  This gives a trivialisation of the local system $\det$.   

Let $f = f_{1} \otimes \ldots f_{n} \in \A^{\otimes n}$.
 
For each $e$, define
$$
\omega_{e} =    K_{l(e)} + \d l(e)  L_{l(e)} \in\br{ \Omega^{\ast}(\Met(\gamma)) \otimes \A^{\otimes 2}}
$$
if $e$ is not an incoming external edge, and
$$
\omega_{e} = (-1)^{p(\A) \abs{f}}(K_{l(e)} + \d l(e)  L_{l(e)}) f_{i} \in \br{ \Omega^{\ast}(\Met(\gamma)) \otimes \A^{\otimes 2}}
$$
if $e$ is an incoming external edge, attached to the incoming external vertex $i$.   The bar indicates that  $\omega_{e}$ is an element of a distributional completion of $ \Omega^{\ast}(\Met(\gamma)) \otimes \A^{\otimes 2}$.

If $e$ is an edge, let $H(e)$ be the two-element set of half-edges of $e$.  Since $e$ is oriented, $H(e)$ is equipped with an isomorphism to $\{0,1\}$, and thus we can think of $\omega_{e}$ as an element 
$$
\omega_{e} \in \br{ \Omega^{\ast}(\Met(\gamma)) \otimes \A^{\otimes H(e)}}
$$
Changing the orientation of $e$ changes the sign of $\omega_{e}$ by $(-1)^{p(\A)}$.

Define
$$
\tilde K_{\gamma} (f) =  \otimes_{e \in E(\gamma)} \omega_{e} \in  \br {\Omega^{\ast}(\Met(\gamma)) \otimes \A^{\otimes H(\gamma)} } $$
where $H(\gamma)$ is the set of half-edges of $\gamma$.  This is an element of the distributional completion of  $\Omega^{\ast}(\Met(\gamma)) \otimes \A^{\otimes H(\gamma)}$. It is singular, because the heat kernel acquires singularities at $t = 0$.

If $v \in V(\gamma)$ is a vertex, and $H_{v}(\gamma)$ is the set of half-edges of $\gamma$ at $v$, then multiplication in the cyclic order on $H_{v}(\gamma)$, combined with the trace map $\Tr : \A \to \C$, defines a map
$$
\Tr_{v} : \A^{\otimes H_{v}(\gamma)} \to \C
$$
We can take the tensor product
$$
\otimes_{v \in V(\gamma) \setminus [m] } \Tr_{v} : \A^{\otimes H(\gamma) \setminus [m]} \to \C
$$
where the tensor product is taken in the chosen ordering on the set $V(\gamma) \setminus [m]$. 

This construction results in an element
$$
K_{\gamma}(f) = \otimes_{v \in V(\gamma) \setminus [m] } \Tr_{v} \tilde K_{\gamma}(f) $$
A priori, this is an element of the distributional completion of $\Omega^{\ast}(\Met(\gamma)) \otimes \A^{\otimes m}$.  
\begin{lemma}
$K_{\gamma}(f)$ is non-singular.    So we have a map
$$
K_{\gamma} : \A^{\otimes n}\to \Omega^{\ast}(\Met(\gamma)) \otimes A^{\otimes m}
$$
This map commutes with the differentials.
\end{lemma}
\begin{proof}
The possible singularities  come from the fact that the heat kernel $K_{t}$ becomes a distribution as $t \to 0$. The conditions on the allowable metrics on $\gamma$, namely that no loop or path spanning outgoing vertices has zero length, imply that $\tilde K_{\gamma}(f)$ is non-singular. 

The heat equation implies that each  $\pi_{e}^{\ast}  K_{l(e)} + \d l(e)\pi_{e}^{\ast} L_{l(e)}$ is closed, which implies that $K_{\gamma}$ commutes with differentials. 

Note that although $K_{\gamma}$ as above is only defined on decomposable elements of $\A^{\otimes n}$, the definition extends to all elements without difficulty. 
\end{proof}

Now if we change the orientation on an edge of $\gamma$, then $ K_{\gamma}$ changes sign by $(-1)^{p(\A)}$.  Similarly if we change the ordering on the sets $E(\gamma), V(\gamma)\setminus [m]$, $K_{\gamma}$  changes sign by the sign of the permutation to the power of $p(\A)$.

Thus we find that  $K_{\gamma}(f)$ defines an element of  $(\Omega^{\ast}(\Met(\gamma),\det^{p(\A)}) \otimes \A^{\otimes m}$, independently of the choice of the ordering on $V(\gamma)$ and $E(\gamma)$.  Here $\det$ is the rank one $\Z/2$ graded local system on $\Met(\gamma)$ defined earlier.   

It is easy to check that the parity of $K_{\gamma}(f)$ is the same as that as $f$, i.e. the map $K_{\gamma}$ is an even map.  The $\Z/2$ grading in the local system $\det$ cancels out the odd parity introduced from the fact that when $p = 1$, the heat kernel $K_{t}$ is odd.

Now $K_{\gamma}$ is clearly $\Aut(\gamma)$ equivariant.  If $f \in \A^{\otimes n}$, we want to show that $K_{\gamma}(f)$ for varying $\gamma$ defines an element of the space of twisted cellular differential forms on the moduli space of graphs, tensored with $\A^{\otimes m}$.  This is guaranteed by the following lemma.
\begin{lemma}
For all $f \in \A^{\otimes n}$,
$$
K_{\gamma}(f)|_{\Met(\gamma/e)}  = K_{\gamma/e}(f) 
$$
\end{lemma}
\begin{proof}
When $t \to 0$, $K_{t}$ becomes the delta distribution along the diagonal.  Now $\Met(\gamma/e) \into \Met(\gamma)$ is the subspace of metrics which give edge $e$ length zero.  Adding on an edge with the delta distribution on it has no effect. 
\end{proof}
Thus, $K_{\gamma}$ for varying $\gamma$ defines a chain map
$$
K : \A^{\otimes n} \to \Omega^{\ast}_{cell}(\Gamma(n,m),\det^{p(\A)}) {\otimes} \A^{\otimes m}
$$
Next we need to check how this behaves with respect to gluing. Let $\gamma_{1},\gamma_{2}$ be graphs, again without metrics, such that $\gamma_{1}$ has $n$ incoming and $m$ outgoing, $\gamma_{2}$ has $m$ incoming and $l$ outgoing edges. Then there is a map
$$
\Met(\gamma_{2}) \times \Met(\gamma_{1}) \to \Met(\gamma_{2}\circ \gamma_{1})
$$
\begin{lemma}
For each $f \in \A^{\otimes n}$, the pull back under this map of $K_{\gamma_{2}\circ \gamma_{1}}(f)$ is 
$$K_{\gamma_{2}} ( K_{\gamma_{1}} (f)) \in \Omega^{\ast}(\Met(\gamma_{2})\times \Met(\gamma_{1}), \det^{p(\A)}) {\otimes} A_{M^{l}} $$
\end{lemma}

Similarly, under the disjoint union map,  $K_{\gamma}$ behaves well.  
Thus we have proved the following theorem.
\begin{theorem}
There is a symmetric monoidal functor from $C_{\ast}(\Gamma,\det^{p(\A)})$ to the category $\Comp$ of chain complexes over $\C$, sending $n \mapsto \A^{\otimes n}$, and a chain $\alpha \in C_{\ast}(\Gamma(n,m),\det^{p(\A)})$ to the map
\begin{align*}
\A^{\otimes n}&\to A^{\otimes m}\\
f &\mapsto \int_{\alpha} K(f) 
\end{align*}
\end{theorem}

\subsubsection{Extension to infinite length edges}
Let $\br \Gamma(n,m)$ be the partial compactification of $\Gamma(n,m)$ obtained by allowing the edges to have infinite length.
\begin{lemma}
The operation $K$ extends to a map
$$K : \A^{\otimes n} \to \Omega^{\ast}_{cell}(\br\Gamma(n,m),\det^{p(\A)}) {\otimes} \A^{\otimes m}$$
\label{lemma_infinite_edge}
\end{lemma}
\begin{proof}
The operator $H$ is a self-adjoint elliptic operator (self-adjoint for the Hermitian pairing on $\A$).  Therefore, in the $L^{2}$ completion of $\A$, there is a basis of eigenvectors $e_{i}$, with corresponding real eigenvalues $\lambda_{i} \ge 0$.  Therefore, the limit as $t \to \infty$ of $e^{-t H}$ is the $L^{2}$ projection onto $\Ker H$.   It follows that the $t \to \infty$ limit of the heat kernel $K_{t}$ is the kernel for projection onto $\Ker H$.

Also $Q^{\dag} t^{-k} e^{-t H} \to 0$  as $t \to \infty$,  for all $k \ge 0$, as $Q^{\dag}$ is zero on harmonic elements of $\A$.   Therefore the kernel $t^{-k} L_{t}$, which represents the operator $Q^{\dag} t^{-k} e^{-t H}$, converges to $0$ as $t \to \infty$.  This implies that $\d t L_{t}$ is zero at $t = \infty$.  

\end{proof}

\section{Topological conformal field theories}
Let me discuss how this construction fits into the more theoretical framework of \cite{Cos04}.   To be self contained, I will start by briefly recalling the definitions of open and closed TCFT.   

Let $\mc M_{\Open}(r,s)$ be the moduli space of Riemann surfaces with $r$ incoming and $s$ outgoing open boundaries.  These surfaces are possibly disconnected.  There are gluing maps $\mc M_{\Open}(r,s) \times \mc M_{\Open}(q,r) \to \mc M_{\Open}(q,s)$, which make these spaces into the morphisms of a topological category, whose objects are $\Z_{\ge 0}$.  This category does not have strict units; however, the disc with one incoming open and one outgoing open boundary is an approximate unit.  We can visualise this disc as a strip with an open boundary on either end.  In order to make $\mc {M}_{\Open}(r,s)$ into a unital category, we allow infinitely short strips.  

We impose the stability condition, that no connected component of a surface in $\mc M_{\Open}(r,s)$ can be a disc with $< 2$ open boundaries.    By moduli space I mean coarse moduli space.   

Let $\mc M_{\Closed}(r,s)$ be the moduli space of Riemann surfaces with $r$ incoming and $s$ outgoing closed boundaries.  As above, these surfaces may be disconnected.  We require that every connected component has at least one incoming boundary. We disallow surfaces which have a connected component which is a disc. As above, these spaces form the morphisms in a topological category, whose objects are $\Z_{\ge 0}$.  This category as defined is not strictly unital; in order to make it strictly unital, we allow infinitely thin annuli.  
 
\begin{remark}
Here we are considering the non-unital version of the theory, as described in section 1.6 of \cite{Cos04}.   The non-unital version of the theory imposes stronger conditions on the surfaces allowed; the unital version allows discs with $\le 1$ open boundary, and discs with one closed boundary.    The non-unital version corresponds of course to non-unital cyclic $A_{\infty}$ categories.

One can show (with a small amount of technical fiddling) that if $\A$ is a unital algebra, as it always is in practise, then the open TCFT we construct is also unital. 
\end{remark}

To these topological symmetric monoidal categories we can associate differential graded symmetric monoidal categories.  Let $C_{\ast}$ denote normalised singular simplicial chains, and let $\Open$ be the category whose objects are $\Z_{\ge 0}$ and whose morphisms are 
$$
\Open(r,s) = C_{\ast}(\mc M_{\Open}(r,s))
$$
We similarly define a dg symmetric monoidal category $\Closed$.
\begin{definition}
A (non-unital) open TCFT (over $\C$) is a dg symmetric monoidal functor $\Open \to \Comp$ to chain complexes over $\C$.  A (non-unital) closed TCFT is a dg symmetric monoidal functor $\Closed \to \Comp$. 
\end{definition}
To deal with odd dimensional theories, for example coming from odd-dimensional Calabi-Yaus, we need to use coefficients in a certain local system $\det$ on these moduli spaces.  Let $\Open^{p},\Closed^{p}$ denote chains with coefficients in $\det^{p}$.

\begin{proposition}
The topological symmetric monoidal category $\Gamma$ is rationally weakly equivalent to the category $\mc M_{\Open}$.

Therefore $C_{\ast}(\Gamma)$ is quasi-isomorphic to $\Open$. More generally $C_{\ast}(\Gamma,\det^{p})$ is quasi-isomorphic to $\Open^{p}$.
\end{proposition}
\begin{proof}

The fact that the space $\Gamma(m,n)$ is homotopy equivalent to $\mc M_{\Open}(m,n)$ follows from the standard ribbon graph decomposition of moduli spaces.   One can  check  that this homotopy equivalence can be made compatible with gluing.  The easiest way to see this is to use the conformal approach to the ribbon graph decomposition, as opposed to the hyperbolic approach.   For each metrised ribbon graph, the associated surface is obtained as follows.   Replace each edge $e$ by $[-1,1]\times [0,l(e)]$.  These are glued together at each vertex, to yield a surface with a flat metric with singularities at the vertices.  This surface has a well-defined conformal structure, and gives a surface with boundary.  This map defines a weak rational homotopy equivalence between the space of metrised ribbon graphs and the space of surfaces with boundary, compatible with gluing. 
\end{proof}

\subsection{Equivalence of two open TCFTs}

We constructed above a symmetric monoidal functor $C_{\ast}(\Gamma,\det^{\otimes p(\A)} \to \Comp$. It follows that for a Calabi-Yau elliptic space $(M,\A)$ with a metric, we get an open TCFT.

In \cite{Cos04} I showed that open TCFTs are the same as a kind of non-unital cyclic (or Calabi-Yau) $A_{\infty}$ category.  One can apply the explicit form of the homological perturbation lemma to the algebra $\A$, as in \cite{Mer99,KonSoi01,Mar04}, to yield a cyclic $A_{\infty}$ structure on $H^{\ast}(\A)$.  The fact that this $A_{\infty}$ algebra is cyclic was proved by Kajiura \cite{Kaj03}  in his thesis.   The condition needed is that the homotopy  between the projector and the identity is skew self adjoint.  (This was explained to me by Ezra Getzler).  To this cyclic $A_{\infty}$ algebra is associated an open TCFT.

\begin{proposition}
These two open TCFTs are quasi-isomorphic.  
\end{proposition}
\begin{proof}[Sketch of proof]
It suffices to show that the two associated cyclic $A_{\infty}$ algebras are quasi-isomorphic.   One is obtained by applying the homological perturbation lemma to $\A$, the other by considering the tree level part of our construction.

We can compactify the moduli space of ribbon trees by allowing infinite length edges.    The $t \to \infty$ limit of the heat kernel is the operator of projection onto harmonic elements of $\A$, i.e.\ the cohomology of $\A$.  As we saw in lemma \ref{lemma_infinite_edge}, the forms extend to these compactifications.

Consider the chain given by all ribbon trees with $n$ incoming vertices, one outgoing vertices, with all external edges of infinite length. When this is oriented correctly, sending $m_{n}$ to this chain gives a quasi-isomorphism between the $A_{\infty}$ operad and the chain operad associated to the moduli space of rooted ribbon trees.   

The integral of our form over this cycle is an operator from $\mc H^{\otimes n} \to \mc H$, where $\mc H$ is the space of harmonic elements of $\A$.  This is because $t \to \infty$ limit of the heat kernel is the operator of projection onto harmonic elements of $\A$.

This integral yields  precisely the formula for $m_{n}$ given by the explicit form of the homological perturbation lemma.  The point is that the integral 
$$\int_{t = 0}^{\infty} \d t  K_{t}$$
gives the Green's operator.  The integral
$$
\int _{t = 0}^{\infty} \d t L_{t}  = - \int _{t = 0}^{\infty} \d t Q^{\dag}  K_{t}  = - Q^{\dag }G
$$
gives the homotopy between the identity map and this projector. The projector and the homotopy are the data needed to construct the explicit version of the homological perturbation lemma; this formula is given by a sum over trees, where at each internal edge we place the homotopy, and each external edge the projector.  This is precisely the same as what we get by integrating our form.

\end{proof}

It follows that the class in $H^\ast(\mc M_g^h)$ associated to a Calabi-Yau elliptic space $(M,\A)$ are the same as those given by Kontsevich's construction \cite{Kon94} applied to cyclic $A_\infty$ algebra $H^\ast(\A)$, when $p(\A) = 0$.    One can see this in more down to earth terms as well.  Let $D_{g,h} \subset \br{\Gamma}_{g,h}$ be the subspace of metrised graphs such that when we cut along all the infinite length edges, what is left is a disjoint union of trees. The inclusion $D_{g,h} \subset \br{\Gamma}_{g,h}$ is a weak homotopy equivalence; this follows from the results of \cite{Cos06}, but it is also easy to prove this directly.

If $\gamma \in D_{g,h}$, let $C(\gamma)$ be the graph obtained by collapsing all edges of $\gamma$ which are not infinite. Note that  $C(\gamma)$ is obtained from $\gamma$ by collapsing a disjoint union of trees.   $D_{g,h}$ is a compact orbi-cell complex; metrised graphs $\gamma,\gamma'$ are in the same cell if $C(\gamma) = C(\gamma')$.  Thus, the cells of $D_{g,h}$ are labelled by ribbon graphs.

If $\gamma$ is a ribbon graph let 
$$
D_\gamma = \{ \gamma' \in D_{g,h} \text{ with an isomorphism }  C(\gamma') \iso \gamma  \}
$$
be the corresponding cell.
Then the integral
$$
\int_{D_\gamma} K_{g,h}
$$
converges.  These integrals define a cellular cocycle on $D_{g,h}$.   The complex of such cellular cochains is just the ribbon graph complex studied in \cite{Kon94}. This cocycle is the same as that obtained by applying Kontsevich's construction to $A_\infty$ algebra $H^\ast(\A)$, with the $A_\infty$ structure given by the explicit form of the homological perturbation lemma.  

\subsection{Constructing closed TCFTs from open}

The results of \cite{Cos04} construct for each open TCFT a closed TCFT. A different construction of closed TCFTs, staring with a cyclic $A_\infty$ algebra, was previously given  by Kontsevich and Soibelman \cite{Kon03,Kon04, KonSoi06}.   Recent work of Tradler and Zeinalian \cite{TraZei06} and Kaufmann \cite{Kau06} give new approaches to this problem. 

I will discuss the approach taken in \cite{Cos04}, which is more abstract than the other constructions.  Let  $F_{M} : \Closed \to \Comp$ be the closed TCFT associated to an open TCFT.    The results of \cite{Cos04} show that 
$$
H_{\ast}(F_{M}(n)) = HH_{\ast}( \A )^{\otimes n}
$$
Here $HH_{\ast}$ can be defined using either the completed projective or ordinary tensor product. We get the same answer, because $\A$ has finite dimensional homology with respect to the operator $Q$.    One technical point is that as we are working with non-unital TCFTs, as in \cite{Cos04} section 1.6, by Hochschild homology we mean the natural Hochschild homology for non-unital algebras.  This is constructed by formally adjoining a unit and then working relative to the subalgebra given by the unit.

The construction of the closed TCFT associated to an open one in \cite{Cos04} is abstract; we define a homotopy-universal open closed theory, and restrict to the closed part of it.    For the reader not familiar with \cite{Cos04}, here is a description of this universal procedure.  
\begin{figure}

\includegraphics[ bb = 1 30 200 200]{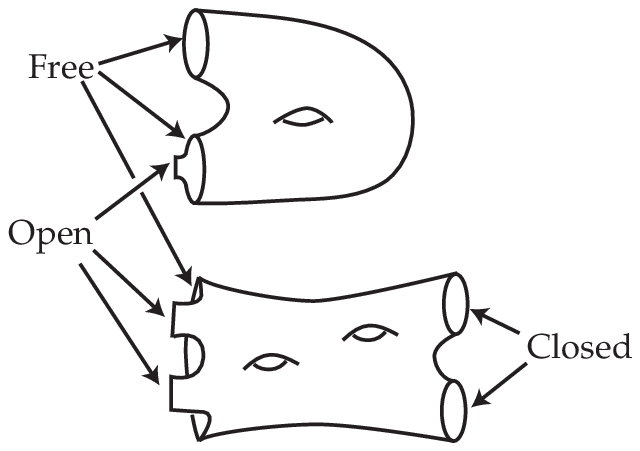}
\caption{ \label{figure surface} }
\end{figure}

Consider a Riemann surface with open and closed  boundaries,  as in figure \ref{figure surface}.   The open boundaries are a disjoint set of closed intervals in $\partial \Sigma$, and the closed boundaries are a disjoint set of parameterised circles in $\partial \Sigma$.  The free boundary of $\Sigma$ is the remainder of the boundary of $\Sigma$.  

Let us suppose, as in the picture, that  the open boundary intervals are all incoming, and the closed boundary circles are all outgoing.  Then one can think of such a surface as describing some open strings joining up to become closed strings.  The initial open strings are given by the open boundary components of the surface.   The free boundary can be thought of as the open string endpoints; this is the part of $\Sigma$ that is constrained to lie on a D-brane.  Note, however, that some boundary components of $\Sigma$ may consist entirely of free boundary.  

Let us denote the moduli space of such surfaces, with $m$ incoming open and $n$ outgoing closed boundaries, by $\mc M_{\OC}(m,n)$.  We require that each connected component of the surfaces have at least one free boundary, and we disallow surfaces which have a connected component which an annulus with one closed boundary and no open boundaries. Let us write
$$
\OC(m,n) = C_{\ast}(\mc M_{\OC}(m,n))
$$

Now suppose we have an open TCFT, given by a symmetric monoidal functor $\mc H_{O} : \Open \to \Comp$.  For each $n \ge 0$, we have a complex of open states $\mc H_{O}(n)$.  This is the complex associated to $n$ open strings.  We can assume, for simplicity,  that $\mc H_{O}(n) = \mc H_{O}(1)^{\otimes n}$ for each $n \ge 0$.  

We can think of $\OC$ as being acted on by the category $\Open$.  Gluing open boundaries of Riemann surfaces together induces maps $\OC(s,n) \otimes \Open(r,s) \to \OC(r,n)$.    Similarly, gluing closed boundaries gives us maps $\Closed(n,m) \otimes \OC(s,n)\to \OC(s,m)$.  This means that we can think of $\OC$ as defining a $\Closed-\Open$ bimodule. 

Then we define
$$
\mc H_{C}(n) = \OC(-,n) \otimes_{\Open}^{\mbb L} \mc H_{O}
$$
The bimodule structure on $\OC$ implies that these spaces carry a natural action of  $\Closed$.

This is a tensor product of co- and contravariant functors from the category $\Open$.  However such tensor products are defined in exactly the same way as tensor products of right and left modules over an algebra.    

Intuitively, we can think of this construction as follows.  The space $\mc H_{C}(n)$ is associated to $n$ closed strings.  Suppose we have a Riemann surface with $m$ incoming open and $n$ outgoing closed boundaries, together with an element of $\mc H_{O}(m)$, the space associated to $m$ open strings.  Then we get an element of $\mc H_{C}(n)$.  More generally, singular chains in the moduli space of such surfaces, labelled in the same way by elements of $\mc H_{O}(m)$,  yield elements of $\mc H_{C}(n)$ in the same way.  Finally we have to impose, in a ``derived'' sense, some natural relations expressing compatibility with gluing surfaces with only open boundary intervals.

One of the results of \cite{Cos04} is that when we do this, the homology of the space $\mc H_{C}(n)$ is the $n$'th tensor power of the Hochschild homology of the $A_{\infty}$ algebra associated to the open theory.

 \section{Chern-Simon type gauge theories}
\label{section cs}

\subsection{Functional integrals over complex vector spaces}

Let $V$ be a finite dimensional complex $\Z/2$ graded vector space, with a non-degenerate, symmetric, complex linear inner product $\ip{\quad}$.  Let $f$ 
be an even (in the $\Z/2$ graded sense) polynomial function on $V$. 

Let $V_\R \subset V$ be a real subspace, such that $V_\R \otimes_\R \C = V$, and such that the inner product is real valued and positive definite on $V_\R$.  Fix a real number $\lambda < 0$. Let $\d \mu$ be any Lebesgue measure on $V_\R$. 
Let
$$
Z (\lambda) =  \int_{x \in V_\R} e^{  \frac{1}{2} \ip{x,x} / \lambda  + f(x) / \lambda  }  \d \mu \left/  \int_{x \in V_\R} e^{  \frac{1}{2} \ip{x,x} / \lambda } \d \mu \right.
$$
Let's suppose that the integrals converge. Wick's lemma allows us to write $Z(\lambda)$ as a sum over graphs, in the standard way, each graph $\gamma$ labelled by $\lambda^{-\chi(\gamma)}$.    We will think of $Z(\lambda)$ as a formal series in $\lambda$.  This series is always defined, whether or not the integral converges. 
\begin{lemma}
The series $Z(\lambda)$ is independent of the choice of real subspace $V_\R \subset V$.
\end{lemma}
\begin{proof}
This follows from Wick's lemma.  Let 
$$P \in V \otimes_\C V$$ 
be minus the inverse tensor to the pairing $\ip{\quad}$ on $V$.    If $v_i$ is a basis for $V$, and
$$
P = \sum p_{ij} v_i \otimes v_j
$$
then
$$
\sum_j p_{ij} \ip{v_j, v_k} = - \delta_{ik}
$$
We can think of $P$ as an element of 
$$V_\R \otimes_\R V_\R \subset V \otimes_\C V$$
Then, $P$ plays the role of the propagator in the Feynman diagram expansion of the integral defining $\ip{f}$.  The interaction at an $i$ valent vertex is defined by the part of $f$ which is homogeneous of degree $i$. 

Since both the propagator $P \in V \otimes_\C V$ and the interactions are independent of the choice of $V_\R$, so is $Z(\lambda)$. 
\end{proof}

In general, we will use the notation
$$
\int_{V} e^{  \frac{1}{2} \ip{x,x} / \lambda+ f/\lambda } \left/   \int_{x \in V} e^{  \frac{1}{2} \ip{x,x} / \lambda } \right. 
$$
to indicate the series obtained by the procedure above : picking any real subspace $V_\R$ on which $\ip{\quad}$ is positive definite, and performing the integral.  Performing the Feynman rules \emph{complex linearly} gives the expansion of this integral in powers of $\lambda$.  By this I mean we think of the propagator $P$ as an element of $V \otimes_\C V$, the interactions as complex linear maps $V^{\otimes_\C k} \to\C$, etc.

If $V$ is infinite dimensional, we will attempt to define such integrals over $V$ by their Feynman diagram expansion, performed complex linearly.

\subsection{ Formal functions }

We want to consider various formal functions on the space $\mc H = \Ker H$ of harmonic elements of $\A$.  Let $R$ be a non-unital nilpotent super-commutative ring.  Such a formal function can be described as something which assigns to each \emph{even} element $a \in  \mc H \otimes R$, an element of $R$, and which is natural with respect to ring homomorphisms $R \to R'$.    In this section, whenever we notation like $Z(a)$, it is to be understood as a formal function of $a\in \mc H$ in this sense.

\subsection{Chern-Simons functional integrals}

Let $\Mat_{N}(\A)$ denote the tensor product algebra $\A \otimes_{\C} \Mat_{N}$, where $\Mat_{N}$ is the algebra of $N \times N$ complex matrices.  The operators $Q,Q^{\dag}$ extend to $\Mat_{N}(\A)$ in the obvious way,
\begin{align*}
Q(a \otimes C ) &= Qa \otimes C \\
Q^{\dag} (a \otimes C ) &= Q^{\dag} a \otimes C
\end{align*}
There is a trace operator 
$$
\Tr : \Mat_{N}(\A) \to \C 
$$
defined by
$$
\Tr (a \otimes C ) = \Tr a  \Tr C
$$
where $\Tr C$ is the usual matrix trace. 

Now assume $p(\A) = 1$.  Consider the Chern-Simons type action on $\Mat_{N}(\A)$, defined by
$$
S(B) =  \frac{1}{2} \Tr B Q B  +  \frac{1}{3} \Tr B^{3}
$$
for $B \in \Mat_N(\A)$.  $S$ is an even functional on the space $\Pi \Mat_{N}(\A)$.

Let $a \in \mc H$, the space of harmonic elements of $\A$. We will try to make sense of the functional integral
$$
Z_{CS}(a,\lambda,N) =  \int_{B \in \Pi \Im Q^{\dag} } e^{S(B + a \otimes \op{Id})/ \lambda} \left/   \int_{B \in \Pi \Im Q^{\dag}} e^{S_k(B)/ \lambda}   \right.
$$
where $S_k(B) = \frac{1}{2} \Tr B Q B$.  The integral is over the super-vector space $\Pi \Im Q^{\dag}\subset \Pi \Mat_{N}(\A)$. This expression is not really meant to be understood as an analytic integral.  Instead it denotes a formal procedure, of expanding over Feynman diagrams, and associating to each Feynman diagram a finite-dimensional integral.    The result (if these integrals converged) would be a formal function of $a \in \mc H$, a formal series in $\lambda$, and an actual function of $N$.  

The condition that $B \in \Im Q^{\dag}$ in the integral is a kind of gauge fixing condition.  There is a direct sum decomposition
$$
\Mat_{N}(\A) = \Im Q^{\dag} \oplus \Ker H \oplus \Im Q
$$
where the subspaces $\Im Q^{\dag}, \Im Q$ are isotropic with respect to the pairing $\Tr(BB')$.  The finite dimensional space $\Ker H$ of harmonic elements is orthogonal to both $\Im Q$ and $\Im Q^{\dag}$.  

In the Batalin-Vilkovisky approach to quantisation (see \cite{Sch92, AleKonSch95}) one integrates over a Lagrangian subspace in the space of fields.  On this Lagrangian, the action functional is non-degenerate.  The choice of Lagrangian is a gauge-fixing condition. The space $\Im Q^{\dag}$ is almost a Lagrangian; it deviates from being a Lagrangian by a finite dimensional subspace.  

The integral we perform is a perturbation expansion around the critical point of the action functional given by $B = 0$.  This critical point is isolated. To see this, observe that  if $\epsilon B$ is a critical point, where $\epsilon$ is a parameter of square zero, then  $Q B = 0$.  However, $\Ker Q \cap \Im Q^{\dag} = 0$.   

Therefore, the integral we consider is the purely ``quantum'' part of Chern-Simons theory.   Non-trivial critical points correspond to deformations of the theory, given by deformations of the operator $Q$ to $Q + B$.  For example, in the case when $M$ is a Calabi-Yau manifold and $\A = \End (E) \otimes _{\cinfty(M)}\Omega^{0,\ast}$, for some holomorphic vector bundle $E$ on $M$, non-trivial critical points correspond to deformations of the vector bundle $E \otimes \C^{N}$.

The propagator for the integral comes from the inverse to the isomorphism 
$$
Q : \Im Q^{\dag} \to \Im Q
$$
This has the integral expression, for $a \otimes C \in \Im Q \subset \Mat_{N}(\A)$,
$$
Q^{-1} (a \otimes C) = \int_{t = 0}^{\infty} Q^{\dag} K_{t} ( a) \d t \otimes C = - \int_{t = 0}^{\infty} L_{t} ( a) \d t  \otimes C
$$
To see this, observe that
\begin{align*}
\int_{t = 0}^{\infty} Q Q ^{\dag} K_{t} (a) \d t &= \int_{t = 0}^{\infty} H K_{t}(a) \d t \\
&= - \int_{t = 0}^{\infty} \frac{\d}{\d t} K_{t}(a) \d t \\
&= a - \lim_{t \to \infty} K_{t}(a)
\end{align*}
The limit $\lim_{t \to \infty} K_{t}$ is the kernel representing orthogonal projection onto $\Ker H \subset \A$.  Since $a  \in \Im Q$, this is zero.

Thus, formally, the propagator is 
$$
P = - \int_{t = 0}^{\infty} L_{t} \otimes \sum X_{ij} \otimes X_{ji}  \in   \A \otimes \A \otimes \Mat_{N} \otimes \Mat_{N} = \Mat_{N}(\A) \otimes \Mat_{N}(\A)
$$
Here $X_{ij} \in \Mat_{N}$ is the matrix whose $ij$ entry is $1$, and all other entries are zero.  This is only a formal expression, as $\int_{t = 0}^{\infty} L_{t}$ is singular. 

Let 
$$P_{\epsilon}  = - \int_{t = \epsilon}^{\infty} L_{t} \otimes \sum X_{ij} \otimes X_{ji}$$
This regularised propagator is non-singular.  The definition of $L_{t}$ implies that 
$
\int_{t = \epsilon}^{\infty} L_{t} 
$
lies in $\Im Q_{1}^{\dag}$.  This expression is also anti-symmetric, which implies that the regularised propagator $P_{\epsilon}$ lies in 
$$
\Im Q^{\dag}\otimes \Im Q^{\dag} \subset \Mat_{N}(\A) \otimes \Mat_{N}(\A)
$$
The usual Feynman rules, which are described very clearly in \cite{Loo93} in the matrix case we are using, show that the regularised Chern-Simons partition function (using the regularised propagator) is
$$
Z_{CS}(a,\lambda,N,\epsilon) = \sum_{\gamma} \lambda^{-\chi(\gamma)} N^{h(\gamma)} \frac{1}{\# \Aut(\gamma)} w(\gamma)
$$
The sum is over trivalent ribbon graphs $\gamma$ possibly with some one-valent external vertices.  The expression $w(\gamma)$ is obtained in the usual way.  On each internal edge we put the propagator $P_\epsilon \in \Mat_{N}(\A)\otimes_\C \Mat_{N}(\A)$.  On each of the $n$ external edges, we put $a$ on the end which connects to the main part of the graph, and ignore the other end. This yields an element of 
$$
\Mat_{N}(\A)^{\otimes_{\C} H(\gamma) \setminus [n] }
$$
(we remove the $n$ half-edges attached to external vertices).  At each internal vertex, we put the linear map
\begin{align*}
\Mat_{N}(\A)^{\otimes_{\C} 3} &\to \C \\
B_{1} \otimes B_{2} \otimes B_{3} &\mapsto \Tr (B_{1} B_{2} B_{3})
\end{align*}
Tensoring these together yields a linear map
$$\Mat_{N}(\A)^{\otimes_{\C} H(\gamma) \setminus [n] } \to \C$$
which we apply to the element constructed from the propagators.

It is easy to see that  if $\gamma$ has $n$ external vertices,
$$
w(\gamma) = \pm \int_{\Met_0^{\epsilon}(\gamma)} K_{\gamma}(a^{\otimes n})
$$
where $\Met_0^{\epsilon}(\gamma)$ is the space of metrics on $\gamma$ in which no internal edge has length $< \epsilon$, and every external edge has length $0$.  $K_{\gamma}(a^{\otimes n})$ is the form on this space with coefficients in the local system $\det$ constructed earlier. The expression  $\int_{\Met_0^{\epsilon}(\gamma)} K_{\gamma}(a^{\otimes n})$ is to be understood as a formal function of $a \in \mc H$ in the sense described above. 

\begin{proposition}
There is a unique isomorphism $\det \iso \Or$ of local systems on $\Gamma_{g,h,n}$ such that sign in this formula is always $+$ (for all CY elliptic spaces of odd parity).
\end{proposition}
\begin{proof}[Sketch of proof]
It's clear that the sign $\pm$ above doesn't depend on the particular CY elliptic space we are using.   We'll construct an orientation with the correct sign for one particular theory, which implies the sign is correct for all theories.    I'll only sketch the proof for $\Gamma_{g,h}$, the proof is similar when we use $\Gamma_{g,h,n}$.

We have already seen (lemma \ref{det_or}) that $\det$ and $\Or$ are isomorphic.    Therefore to give an isomorphism $\det \iso \Or$ is suffices to give it away from codimension $2$ strata.     That is, it suffices to give, for each trivalent ribbon graph $\gamma$, an isomorphism between $\det$ and $\Or$ on $\Met(\gamma)$, which is $\op{Aut}(\gamma)$ equivariant, and which satisfies some compatibility condition on ribbon graphs with a single $4$ valent vertex.    

To give the isomorphism for trivalent graphs, it suffices to give a real section of $\Omega_{cell}^{6g-6+3h+n}(\Gamma_{g,h,n}, \det)$,  which never vanishes on the space $\Met(\gamma)$ for any trivalent ribbon graph $\gamma$.

To construct such a section we write down a simple \emph{real} Calabi-Yau elliptic space, of odd parity, such that the top degree component of the form on moduli space never vanishes.

One such theory is given by taking $M$ to be a point, and $\A = \R[x]/ (x^3 )$, where $x$ is of degree $1$. The trace is given by $\Tr x^3 = 1$, and the differential is $Qx = -x^2$, $Qx^k = 0$ unless $k = 1$.   The operator $Q^\dag$ is $Q^\dag x^2 = -x$, and $Q^\dag x^k = 0$ unless $k = 2$.    The Hamiltonian $H$ is $H x = x$, $H x^2 = x^2$, $H 1 = Hx^3 = 0$.   The resulting Chern-Simons type action on $\Pi \A^1 = \Pi \Im Q^\dag$ is $-\frac{1}{2} x^2 + \frac{1}{3} x^3$. 

Let us take the top-degree (in fact the only non-zero piece) of the form associated to this theory, in $\Omega_{cell}^{6g-6+3h}(\Gamma_{g,h}, \det)$.  It's easy to calculate that this does not vanish on $\Met(\gamma)$ for any trivalent graph $\gamma$.   

The isomorphism $\det \iso \Or$ on the space of trivalent graphs is chosen so that the integral of this  form over any top-dimensional compact piece is always positive.  

It remains to check the compatibility condition on ribbon graphs with a single $4$ valent vertex.  This is a fairly simple calculation, which I will not perform here.  
\end{proof}

Let $\Gamma_{g,h} \subset \Gamma(0,0)$ be the space of connected metrised ribbon graphs  with $h$ unordered boundary components and Euler characteristic $-2g+2-h$.  More generally, let $\Gamma_{g,h,n} $ be the space of such ribbon graphs which have $n$ ordered external vertices, with the property that the lengths of the edges attached to the external vertices is $0$, and that the external vertices are a non-zero distance apart.  We only consider the case when $2g-2+h + \frac{n}{2} > 0$.

Let $\Gamma^{\epsilon}_{g,h,n} \subset \Gamma_{g,h,n}$ be the subspace of graphs with no interior edges of length $< \epsilon$. For $a \in \A$,  let us denote by $K_{g,h}(a^{\otimes n})$ the restriction of the form $K(a^{\otimes n})$ in $\Gamma(n,0)$ to $\Gamma_{g,h,n}$.

The discussion above shows that 
$$
Z_{CS}(a,\lambda,N,\epsilon) =  \exp \left (\sum_{\substack {g,n \ge 0 h > 0 \\ 2g-2+h + \frac{n}{2} > 0}} \lambda^{2g-2+h} N^{h} \int_{ \Gamma^{\epsilon}_{g,h,n}} K_{g,h} (a^{\otimes n}) / n! \right)
$$
This expression may not converge as $\epsilon \to 0$.  This is an equality of formal functions of $a \in \mc H$ with values in formal series in $\lambda$, for each value of $N,\epsilon$.

The space $\Gamma_{g,h}$ is homeomorphic to the space $\mc M_g^h$ of Riemann surfaces of genus $g$  with $h$ unordered, unparameterised boundary components.  
We can formally set $\epsilon = 0$ and write down the equality of non-convergent expressions
$$
Z_{CS}(a,\lambda,N) =  \exp \left (\sum_{\substack {g,n \ge 0, h > 0 \\  2g-2+h + \frac{n}{2} > 0}} \lambda^{2g-2+h } N^{h} \int_{\mc M_g^{h,n}} K_{g,h}(a^{\otimes n}) / n! \right)
$$
which we interpret as an identity between an open-string and a Chern-Simons partition function.

\section{Independence of choice of Hermitian metric}
The whole construction is independent, up to homotopy, of deformations of the Hermitian metric on $\A$. The metric is encoded in a complex anti-linear Hodge star operator $\ast : \A \to \A$, such that $\ast^{2} = \pm 1$.

We will show that if we perturb $\ast$, then the deformation of the theory extends to one over a contractible commutative differential algebra, so that we get a homotopic theory.   This implies that when we change $\ast$, all of our forms on moduli space change by the addition of an exact form.  

So suppose $s$ is an even infinitesimal parameter, $s^{2}= 0$, and we perturb $\ast$ to 
$$\ast_{s} = \ast + s  \ast'.$$  Let us suppose this perturbed operator still satisfies all the necessary conditions.

Let $\epsilon$ be an odd parameter, and consider the contractible algebra $ R =\C[s,\epsilon]/(s^{2}, \epsilon s )$, with differential $Q s = \epsilon$.   We want to show that the deformation of the theory over the $\C[s]/s^{2}$ extends to one over $R$.  

All we do is set
$$
\ast_{s,\epsilon} = \ast_{s} = \ast + s \ast'
$$
Then, by assumption, this satisfies all the necessary conditions, so the theory extends.  We have
\begin{align*}
Q^{\dag}_{s,\epsilon} &= Q^{\dag} + s \delta Q^{\dag}  = Q^{\dag} \pm s (  (\ast') Q \ast + \ast Q (\ast'))\\
Q_{s,\epsilon} &= Q + \epsilon \frac{\partial}{\partial s}\\
H_{s,\epsilon} &= [Q,Q^{\dag}_{s,\epsilon}] = H + s [Q,\delta Q^{\dag}] + \epsilon \delta Q^{\dag}
\end{align*}

Note that the operator $H_{s,\epsilon}$ admits a heat kernel if $H$ did, because $H_{s,\epsilon}$ is a perturbation of $H$ by nilpotent operators.

\section{Acknowledgements}
I'd like to thank Maxim Kontsevich for some very helpful correspondence, and in particular for explaining to me his constructions of elliptic spaces; and Mohammed Abouzaid, Paul Seidel and Boris Tsygan for many interesting conversations related to this work.

\section{Appendix}

\subsection{ Identities satisfied by the heat kernel}
Here we prove 
\begin{proposition}[Proposition \ref{prop heat kernel identities}]
The identities
\begin{align}
Q_{x} K_{t}(x,y) +  Q_{y} K_{t}(x,y) &= 0 \\
Q_{x}^{\dag} K_{t}(x,y) &=  Q_{y}^{\dag} K_{t}(x,y)  \\
H_{x}K_{t}(x,y) &= H_{y} K_{t}(x,y)  \\
K_{t}(x,y) &= (-1)^{p}K_{t}(y,x) \\
L_{t}(x,y) &= (-1)^{p} L_{t}(y,x) \\
(Q_{x}+Q_{y})L_{t} (x,y) &= -H_{x} K_{t}(x,y) = \frac{\d}{\d t} K_{t}(x,y)
\end{align}
are satisfied.
\end{proposition}
\begin{proof}
Let $\A^{0}$, $\A^{1}$ refer to the odd and even parts of $\A$.  I will prove these identities in the case when $p = 1$; the even case is easier.

Let us write 
$$
K_{t}(x,y) = K_{t}^{1,0}(x,y) + K_{t}^{0,1}(x,y) \in \A^{1}\otimes  \A^{0} \oplus \A^{0}\otimes \A^{1}
$$
Let $f \in \A$ be a function.
$$
(K_{t} f ) (x) =  (-1)^{\abs{f}}\Tr_{y}(  K_{t}^{1,0}(x,y) + K_{t}^{0,1}(x,y) ) f(y)
$$
Observe that 
\begin{align*}
\Tr_{y} Q_{y} K^{1,0}(x,y) f (y) &= \Tr_{y} K^{1,0}(x,y) Q_{y} f(y) \\
\Tr_{y} Q_{y} K^{0,1}(x,y) f (y) &= \Tr_{y} K^{0,1}(x,y) Q_{y} f(y) 
\end{align*}
These identities use the fact that $Q$ is skew self adjoint. 

This shows that 
$$
(-1)^{\abs f}\Tr_{y} Q_{y} K_{t}(x,y) f(y) = (-1)^{\abs f} \Tr_{y} K_{t}(x,y) Q_{y} f(y) = - ( K_{t}  Q f )(x)
$$
It is clear that 
$$( Q K_{t} f ) (x) = \Tr_{y} Q_{x} K_{t}(x,y) f(y)$$
Since $[K_{t},Q] = 0$ we see that 
$$
\Tr_{y} (Q_{x} + Q_{y} ) K_{t} (x,y) f(y) = 0
$$
as desired. 

Similar arguments prove the second, third, and sixth identities. The different signs arises as $Q^{\dag}, H$ are self adjoint instead of skew adjoint. 

We will again prove the identity $K_{t}(x,y) = (-1)^{p}K_{t}(y,x)$ in the case $p = 1$. Let $K^{\ast}_{t}$ be the operator associated to the kernel $K_{t}(y,x)$.     In more formal terms, let $\sigma : \A \otimes \A \to \A \otimes \A$ be the interchanging operator; then $K_{t}(y,x) = \sigma K_{t}(x,y)$.  It is easy to see that for all $a \in \A \otimes \A$, 
$$
\Tr_{\A \otimes \A} a = -\Tr_{\A \otimes \A} \sigma a
$$
i.e. that changing the order of integration affects things by a sign.
Let $f,g \in \A$, with $f \in A^{0}$ and $g \in A^{1}$.  It follows that 
\begin{align*}
\ip{f, K_{t}^{\ast} g } &=  -\Tr_{\A \otimes \A} f(x) K_{t}(y,x) g(y) \\
&= \Tr_{\A \otimes \A}f(y) K_{t}(x,y) g(x)\\
&= - \Tr_{\A \otimes \A} g(x) K_{t}(x,y) f(y) \\
&= - \ip{g, K_{t}  f}
\end{align*}
As $K_{t}$ is self adjoint and even, we see that $K_{t}^{\ast} = - K_{t}$ as desired. 

The fact that $L_{t}(x,y) = (-1)^{p} L_{t}(y,x)$ follows from this and the second identity.

\end{proof}

\newcommand{\etalchar}[1]{$^{#1}$}
\def\cprime{$'$}

\end{document}